\documentclass[a4paper,10pt,intlimits,sumlimits]{amsart}
\usepackage[utf8]{inputenc}

\usepackage{amsmath}
\usepackage{amsthm}
\usepackage{amstext}
\usepackage{amssymb}
\usepackage{amsfonts}
\usepackage{mathrsfs}
\usepackage{mathtools}

\usepackage{esint}

\usepackage{mathbbol}
\usepackage{bm}
\usepackage{hyperref}

\usepackage{tikz}
\usepackage{pgfplots}
\usetikzlibrary{patterns,calc}

\usepackage{float}
\usepackage{graphicx}
\usepackage{wrapfig}

\usepackage[disable]{todonotes}

\usepackage[citestyle=alphabetic,bibstyle=alphabetic,maxalphanames=5,backend=bibtex]{biblatex}

\newcommand\num{\stepcounter{equation}\tag{\theequation}}
\def\mcl_#1{\limits_{\mathclap{#1}}}

\newcommand{\mc}{\mathcal}

\newcommand{\mb}[1]{\mathbb{#1}}
\newcommand{\mf}[1]{\mathfrak{#1}}

\newcommand{\R}{\mb{R}}
\newcommand{\C}{\mb{C}}
\newcommand{\N}{\mb{N}}
\newcommand{\Z}{\mb{Z}}

\newcommand{\dd}{\mathrm{d}}
\newcommand{\sthat}{\,\colon\,}

\newcommand{\Fourier}{\mc{F}}
\newcommand{\FT}{\widehat}
\newcommand{\1}{\mb{1}}

\renewcommand{\phi}{\varphi}
\renewcommand{\epsilon}{\varepsilon}

\renewcommand{\tilde}{\widetilde}

\newcommand{\closure}[1]{\overline{#1}}

\DeclareMathOperator{\dist}{dist}
\DeclareMathOperator{\spt}{spt}

\DeclareCiteCommand{\citejournal}
  {\usebibmacro{prenote}}
  {\usebibmacro{citeindex}%
    \usebibmacro{journal}}
  {\multicitedelim}
  {\usebibmacro{postnote}}

\renewcommand{\L}{\hbox{\raisebox{0.06em}-}\kern-0.45emL}
\renewcommand{\a}{\mf{a}}
\renewcommand{\c}{\mf{c}}

\newtheorem{thm}{Theorem}
\newtheorem*{thm*}{Theorem}

\newtheorem{defn}[thm]{Definition}
\newtheorem*{defn*}{Definition}

\newtheorem{prop}[thm]{Proposition}
\newtheorem*{prop*}{Proposition}

\newtheorem{cor}[thm]{Corollary}
\newtheorem*{cor*}{Corollary}

\newtheorem{lem}[thm]{Lemma}
\newtheorem*{lem*}{Lemma}

\newtheorem{rem}[thm]{Remark}
\newtheorem*{rem*}{Remark}

\numberwithin{equation}{section}
\numberwithin{thm}{section}

\allowdisplaybreaks

\begin{document}

\title{Variational Carleson embeddings into the upper $3$-space.} 
\date{\today} 
\author{Gennady Uraltsev}

\begin{abstract}
  In this paper we formulate embedding maps into time-frequency space
  related to the Carleson operator and its variational counterpart. We
  prove bounds for these embedding maps by iterating the outer measure
  theory of \cite{do2015lp}. Introducing iterated outer $L^{p}$
  spaces is a main novelty of this paper.
\end{abstract}

\maketitle

\section{Introduction}

In this paper we consider the Carleson Operator
 \begin{align*}\num  \label{eq:carleson-op}
    \mc C_{c} f (z) := \int_{c(z)}^{+\infty} \FT f (\xi) e^{i\xi z} \dd \xi,
 \end{align*}
 with $c:\R\to \R$ a Borel-measurable stopping function. The
Variational Carleson Operator studied by
 Oberlin et al. in \cite{oberlin2012variation} is given by:
\begin{align*}\num \label{eq:var-carleson-op}
  \mc V^{r}\mc C_{\mf c} f(z)= 
  \left(  \sum_{k\in\Z}\left|\mc C_{\mf c_{k+1}}f(z)-\mc C_{\mf c_{k}}f(z)\right|^{r}
  \right)^{1/r}
\end{align*}
where $\mf c: \Z\times\R \to \R\cup\{+\infty\}$ is a stopping sequence
of Borel-measurable functions such that
$\mf c_{k}(z)\leq\mf c_{k+1}(z)$ for all $z\in\R $ and $k\in\Z$.  The
boundedness on $L^{p}\left(\R\right)$ with $p\in(1,\infty)$ of these
operators, uniformly with respect to the stopping functions $c$ and
$\mf c$, implies the famous Carleson Theorem on the almost everywhere
convergence of the Fourier integral for functions in $L^{p}(\R)$. The
main technique for bounding these operators were first introduced by
Carleson in his paper \cite{carleson1966convergence} on the
convergence of Fourier series for $L^{2}\left([-\pi/2,\,\pi/2)\right)$
functions and is often referred to as time-frequency analysis.

The purpose of this paper is to discuss  embedding maps into
time-frequency space $\mb X=\R\times\R\times \R^{+}$ relevant to
 \eqref{eq:carleson-op} and \eqref{eq:var-carleson-op}. In Theorems
\ref{thm:simple-mass-boundedness}, \ref{thm:energy-boundedness}, and
\ref{thm:max-var-mass-bounds} we show the boundedness properties of
these embedding maps in terms of appropriately defined norms.  Generally
speaking an embedding map is a representation of a function by another
function defined on the symmetry group of the problem at hand. The
appropriate norms for dealing with these embedded functions are the
outer measure $L^{p}$ norms introduced in \cite{do2015lp} in the
context of the Bilinear Hilbert Transform, an operator with the same
symmetries as \eqref{eq:carleson-op} and
\eqref{eq:var-carleson-op}.

Theorem \ref{thm:energy-boundedness} is an extension of the result of
\cite{do2015lp} to $1<p<2$. For our proof we introduce iterated, or
semi-direct product, outer measure $L^{p}$ spaces and incorporate the
idea by Di Plinio and Ou \cite{di2015modulation} of using
multi-frequency Calderón-Zygmund theory from
\cite{nazarov2010calderon}. The embedding Theorems
\ref{thm:simple-mass-boundedness} and \ref{thm:max-var-mass-bounds}
are somewhat dual to \ref{thm:energy-boundedness} for the purpose of
bounding the bilinear form associated to \eqref{eq:carleson-op} and
\eqref{eq:var-carleson-op} respectively.

In \cite{oberlin2012variation} the operator \eqref{eq:var-carleson-op} has
been shown to be bounded for $p\in(1,\infty)$ and $r\in(2,p')$. The
proof in the range $p\in(2,r)$ requires only theorems that make use of
non-iterated outer measure spaces of \cite{do2015lp}.  While initially
introduced only to address the range $p\in(r',2]$, iterated outer
measure spaces surprisingly provide a direct proof in the complete
range $p\in [r,\infty)$, and hereby explain ad-hoc interpolation
techniques used in \cite{oberlin2012variation}.

\begin{wrapfigure}{r}{0pt}
  \begin{tikzpicture}[scale=2.5]

    \def\xmax{1.3} \def\xmin{-0.1} \def\ymax{1.3} \def\ymin{-0.1}
    \draw [->, thick] (0, \ymin) -- coordinate (y axis mid) (0,\ymax);
    \draw [->, thick] (\xmin, 0) -- coordinate (x axis mid) (\xmax,0);
    \node [right] at (1.3,1.1) {\small $\begin{aligned} &p\in(r' , \infty)
        \\ & r\in(2,\infty]\end{aligned}$};

    \fill[color=olive,opacity=0.3] (0,0) -- (0.5,0.5) -- (0,0.5) --
    cycle;
    \fill[color=olive,opacity=0.3] (0.5,0) -- (0.5,0.5) --
    (1,0) -- cycle;

    \fill[color=olive,opacity=0.3] (0,0) -- (0.5,0.5) -- (0.5,0) --
    cycle;
    \fill[color=olive, pattern=dots, pattern color=black] (0,0)
    -- (0.5,0.5) -- (0.5,0) -- cycle;

    \draw[fill=black] (1,0) circle (0.01) node [below] {$1$};
    \draw[fill=black] (0,1) circle (0.01) node [left] {$1$};

    \draw[fill=black] (0.5,0) circle (0.01) node [below]
    {$\frac{1}{2}$}; \draw[fill=black] (0,0.5) circle (0.01) node
    [left] {$\frac{1}{2}$};

    \node[below] at (1.3,0) {$\frac{1}{p}$}; \node[left] at (0,1.3)
    {$\frac{1}{r}$};

    \node [right] at (0.3,1.0) {\small $p\in(r,\infty)$}; \draw
    [->,very thin, dashed](0.3,1.0) to [out=-150, in=150] (0.2,0.4);

    \node [right] at (1.0,0.7) {\small $p\in(2,r)$}; \draw [->, very
    thin,dashed](1.0,0.7) to [out=180,in=30] (0.45,0.3);

    \node [right] at (1.0,0.3) {\small $p\in(r',2)$}; \draw [->,very
    thin,dashed](1.0,0.3) to [out=180,in=30] (0.65,0.2);

    \draw[black,thin] (0,0.5) -- (0.5,0.5) -- (1,0) -- (0,0) -- cycle;


  \end{tikzpicture}
  \caption{Bounds of $\mc V^{r}\mc C_{\mf c}$ on $L^{p}(\R)$.}
  \label{fig:var-carleson-exponents}
\end{wrapfigure}
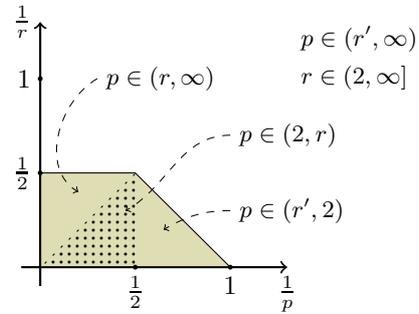

The advantage of reasoning in terms of embedding maps is also attested
by the recent developments in \cite{culiuc2016domination} that prove
sharp weighted bounds for the Bilinear Hilbert Transform using the
embedding from \cite{di2015modulation}. This is done by dominating the
trilinear form associated to the operator by sparse forms following the
approach of Lacey \cite{lacey2015elementary}. In a similar spirit,
sparse domination and weighted boundedness results for
the Variational Carleson Operator are forthcoming in a paper by the Di
Plinio, Do, and the author, that make use of the embedding maps of the
present paper.

We also point out the recent paper \cite{do2016variational} in which
Do, Muscalu, and Thiele use outer measure $L^{p}$ spaces to provide
variational bounds for bilinear  Fourier inversion integrals,
that are bilinear versions of \eqref{eq:var-carleson-op}.

On a historical note, we point out Hunt's extension
\cite{hunt1968convergence} to $L^{p}$ with $p\in(1,\infty)$ of
Carleson's pointwise almost-everywhere convergence result
\cite{carleson1966convergence} for Fourier series of functions on
$L^{2}\left([-\pi/2,\,\pi/2)\right)$. Carleson's and Hunt's results
depend on a fine analysis of the properties of a function on the
torus.  In \cite{fefferman1973pointwise} Fefferman concentrated on
proving the same result by a careful study of the operator
\eqref{eq:carleson-op}. The wave-packet representation for the
operator that is crucial for making use of embedding maps appeared in
\cite{lacey2000proof} that provides a more symmetric approach
encompassing the aforementioned two ideas. This approach inspired both
\cite{oberlin2012variation} and the present paper.

Finally, we emphasize that we formulate an embedding map into the
time-frequency space parameterized by continuous parameters, in the
vein of \cite{do2015lp}. This allows us to avoid model-sum operators
and averaging procedures ubiquitous in other works in time-frequency
analysis. Furthermore, such a formulation proves to be more versatile
and in particular the results of the present paper imply all the
bounds for the discretized model used in \cite{oberlin2012variation}.


\subsection{The Carleson operator}

For simplicity we begin by discussing the Carleson operator
\eqref{eq:carleson-op} that is a specific instance of
\eqref{eq:var-carleson-op} for $r=+\infty$. The operator is given
pointwise by the Fourier multiplier operator associated to the
multiplier $\1_{[c(z),\,+\infty)}(\xi)$ applied to $f$. This can be
expressed in terms of a wavelet frame centered at frequency $c(z)$
using a continuous Littlewood-Paley decomposition:
\begin{align*}\num\label{eq:C-LP}
  \mc C_{c} f(z) = \int_{\R^{+}}\int_{\R}f*\psi_{\eta,t}*
  \psi_{\eta,t}(z)\,\chi\left(t(\eta-c(z))\right) \dd \eta \dd t
\end{align*}
where
\begin{align*}\num\label{eq:psi-eta-t}
  &  \psi_{\eta,t}(z):=t^{-1}e^{i\eta
    z}\,\psi\left(\frac{z}{t}\right)
\end{align*}
with $\psi\in S(\R)$ a suitably normalized, non-negative, even, 
generating wavelet with Fourier transform $\FT \psi$ supported in a
small ball $B_{b}$. We use the notation $B_{r}(x):=(x-r,x+r)$ to denote
a ball of radius $r$  centered at $x$, while if $x=0$ we omit it by
simply writing $B_{r}$. The non-negative cutoff
function $\chi$ satisfies
\begin{align*}\num\label{eq:def:chi-bump}
  & \chi\in C^{\infty}_{c}\left(B_{\epsilon}(d)\right)&& B_{\epsilon}(d)\subset\left(b,\,+\infty\right)&&\int\chi=1.
\end{align*}
Given two functions $f,a\in S(\R)$ set 
\begin{align*}\num\label{eq:F}
 & F(y,\eta,t) := f*\psi_{\eta,t}(y)\\
\num \label{eq:A}
  &A(y,\eta,t) := \int_{\R} a(z)\psi_{\eta,t}(y-z)\,\chi\left(t(\eta-
    c(z))\right) \dd z .
\end{align*}
The arguments of the above functions are points of the time-frequency
space $\mb X=\R\times\R\times\R^{+}$ that parameterizes the defining
symmetries of the class of operators defined by (\ref{eq:carleson-op})
i.e. translation of the function, translation of its Fourier
transform, and dilation.  The outer measure $L^{p}$ spaces allow one
to deal with the overderminancy of the wave-packets.

The wave packet representation \eqref{eq:C-LP} gives the
inequality
\begin{align*}\num \label{eq:duality-wave-packet}
  \left|\int_{\R} C_{c} f(z)\,a(z)\dd z\right| \leq
  \left|\iiint_{\mb X} F(y,\eta,t)A(y,\eta,t)\,
  \dd  y \dd \eta \dd t\right|.
\end{align*}
By duality the bound of the operator \eqref{eq:carleson-op} on $L^{p}(\R)$
follows from bounds on $L^{p}(\R)\times L^{p'}(\R)$ of the bilinear form on the
left hand side of the previous display.

The abstract framework of outer measure $L^{p}$ spaces provides us with the
Hölder type bound
\begin{align*}\num \label{eq:outer-duality}
  \left|\iiint_{\mb X} F(y,\eta,t)A(y,\eta,t)\,\dd \eta
  \dd y \dd t\right|
  \lesssim
  \|F\|_{L^{p}\L^{q}(S_{e})}
  \|A\|_{L^{p'}\L^{q'}(S_{m})}
\end{align*}
with $\frac{1}{p}+\frac{1}{p'}=1$ and
$\frac{1}{q}+\frac{1}{q'}=1$. Appearing on the right are
iterated outer $L^{p}$ quasi-norms that we elaborate on in Section
\ref{sec:outer-measure-def}. 

The embedding maps defined via equations \eqref{eq:A} and
\eqref{eq:F}, that we call ``mass'' and ``energy'' embeddings for
historical reasons (compare with
\cite{lacey2000proof}), satisfy the bounds
\begin{align*}\num\label{eq:sm-bounds}
  &\|A\|_{L^{p'}\L^{q'}(S_{m})}\lesssim\|a\|_{L^{p'}}, \\\num \label{eq:e-bounds}
  &\|F\|_{L^{p}\L^{q}(S_{e})}\lesssim\|f\|_{L^{p}}.
\end{align*}

\begin{thm}[Mass embedding bounds] \label{thm:simple-mass-boundedness}

  For any $p'\in(1,\infty]$, $q'\in(1,\infty]$, and for any function
  $a\in L^{p'}(\R)$ the bounds \eqref{eq:sm-bounds} for the embedding
  \eqref{eq:A} hold with a constant independent of the Borel
  measurable function $c:\R\to\R$.
\end{thm}

\begin{thm}[Energy embedding bounds]\label{thm:energy-boundedness}
  For any $p\in(1,\infty]$, $q\in\left(\max(2;p'),\infty\right]$, and
  for any $f \in L^{p}(\R)$ the bounds \eqref{eq:e-bounds} for the
  embedding \eqref{eq:F} hold. 
\end{thm}

Theorem \ref{thm:simple-mass-boundedness} follows as a corollary of
Theorem \ref{thm:max-var-mass-bounds} below while Theorem
\ref{thm:energy-boundedness} will be proven in Section \ref{sec:e-bounds}.

The boundedness of the Carleson Operator on $L^{p}(\R)$ follow as a
result of the above discussion. Indeed for any $p,p'\in(1,\infty)$ with
$\frac{1}{p}+\frac{1}{p'}=1$ one can find $q,q'\in(1,\infty)$ such
that $\frac{1}{q}+\frac{1}{q'}=1$ and bounds \eqref{eq:sm-bounds} and
\eqref{eq:e-bounds} hold.

We remark that iterated outer measure spaces are used to address the
case $p\in(1,2)$. In Section \ref{sec:e-bounds} we show that if
$p\in(2,\infty)$ a the non-iterated version of outer measure $L^{p}$
spaces are sufficient to prove $L^{p}$ boundedness of
\eqref{eq:carleson-op}.


\subsection{The variational Carleson operator}

The operator \eqref{eq:var-carleson-op}, introduced and studied in
\cite{oberlin2012variation}, is bounded on $L^{p}(\R)$ for
$r\in(2,\infty]$ and $p \in (r',\infty)$. The above paper also shows
that this range is sharp in the sense that that strong $L^{p}$ bounds
do not hold outside this range (see Figure \ref{fig:var-carleson-exponents}).

By duality it is sufficient to prove the bilinear a priori bound
\begin{align*}  \num\label{eq:var-duality-linear}
  \left| \int_{\R} \sum_{k\in\Z} \a_{k}(z)
  \int_{\c_{k}(z)}^{\c_{k+1}(z)} \FT f(\xi) e^{i\xi z}\dd \xi
  \;\dd z \right|
\lesssim \|f\|_{L^{p}}\|\a\|_{L^{p'}(l^{r'})}.
\end{align*}
with a constant independent of the stopping sequence $\c$.  For the above expression to
make sense we require that $f\in S(\R)$ while
$\a\in L^{p'}(l^{r'})$ i.e. $z\mapsto \a(z)=(\a_{k}(z))_{k\in\Z}$ is a
function on $\R$ such that for every $z\in\R$ its  value is the
sequences $\a(z)=(\a_{k}(z))_{k\in\Z} \in l^{r'}(\Z)$. The function
$\a$ is Borel measurable in Bochner sense and
\begin{align*}
  \|\a\|_{L^{p'}(l^{r'})}:=\left(\int_{\R}\|\a(z)\|_{l^{r'}}^{p'}\dd
  z\right)^{1/p'}<\infty.
\end{align*}

Analogously to \eqref{eq:duality-wave-packet}, the left had side of
\eqref{eq:var-duality-linear} admits a wave-packet domination
\begin{align*}\num\label{eq:var-duality-wavepacket}
  \left| \int_{\R} \sum_{k\in\Z}\a_{k}(z)
  \int_{\c_{k}(z)}^{\c_{k+1}(z)} \FT f(\xi) e^{i\xi z}\dd \xi
  \;\dd z \right|\leq \iiint_{\mb X} \left|F(y,\eta,t)\mb
  A(y,\eta,t)\right| \dd y\dd \eta \dd t.
\end{align*}
where the embedding map $\a\mapsto \mb A$ is given by
  \begin{align*}\num\label{eq:def:var-mass-embedding}
    \mb A(y,\eta,t) := \sup_{\Psi} \left|
    \int_{\R}\sum_{k\in \Z}
    \a_{k}(z) \Psi_{y,\eta,t}^{\mf c_{k}(z),\mf c_{k+1}(z)}(z)\dd z\right|.
  \end{align*}
  The supremum above is taken over all possible choices of left or
  right truncated wave packets $\Psi_{y,\eta,t}^{c_{-},c_{+}}$.  A
  left truncated wave packet $\Psi_{y,\eta,t}^{c_{-},c_{+}}$ at
  $(y,\eta,t)\in\mb X$ is a $S(\R)$ function parameterized by
  $c_{-}<c_{+}\in\R\cup\{+\infty\}$. The parameterization satisfies the properties
  below.
  The following three functions of the variable $z$
      \begin{gather*}\num\label{eq:wave-packet-smoothness}
         e^{-i\eta(y+tz)} t\Psi_{y,\eta,t}^{c_{-},c_{+}}(y+tz)
        \\
         t^{-1}\partial_{c_{-}}\Big(e^{-i\eta(y+tz)}
        \Psi_{y,\eta,t}^{c_{-},c_{+}}(y+tz)\Big)
            \\
        t^{-1}\partial_{c_{+}}\Big(e^{-i\eta(y+tz)}
        t\Psi_{y,\eta,t}^{c_{-},c_{+}}(y+tz)\Big)
  \end{gather*} 
  are bounded in $S(\R)$ uniformly for all $(y,\eta,t)\in \mb X$ and
  $c_{-}<c_{+}\in\R$. For some constant $b>0$ the functions
  $\Psi_{y,\eta,t}^{c_{-},c_{+}}$ satisfy
  \begin{align*}\num\label{eq:wave-packet-freq-support}
    &\spt \FT \Psi_{y,\eta,t}^{c_{-},c_{+}} \subset B_{t^{-1}b}(\eta).
  \end{align*}
For some constants $d,d',d''>0$, and $\epsilon>0$ it holds that
  \begin{align*}\num\label{eq:eq:wave-packet-geomtery-vanish}
    &\Psi_{y,\eta,t}^{c_{-},c_{+}}\neq 0
    &&\text{only if } \left\{
       \begin{aligned}
         &t(\eta - c_{-}) \in B_{\epsilon}(d) \\&t(c_{+}
         -\eta)> d'>0
       \end{aligned}
\right.
    \\ \num\label{eq:eq:wave-packet-geomtery-constant}
    &  \Psi_{y,\eta,t}^{c_{-},c_{+}}=
      \Psi_{y,\eta,t}^{c_{-},+\infty}
    &&\text{\phantom{only }if } \phantom{\Bigg\{}t(c_{+}-\eta)> d''>d'>0.
  \end{align*}
The wave packet $\Psi_{y,\eta,t}^{c_{-},c_{+}}$ is right truncated if
$\Psi_{y,-\eta,t}^{-c_{+},-c_{-}}$ is left truncated. 

The main result of this paper is the following bounds for the
embedding \eqref{eq:def:var-mass-embedding} that are analogous to
the bounds \eqref{eq:sm-bounds}.

\begin{thm}[Variational mass
  embedding bounds]\label{thm:max-var-mass-bounds}
  For any $r'\in[1,2)$, $p'\in(1,\infty]$, and $q'\in(r',\infty]$ and
  any function  $\a\in L^{p'}(l^{r'})$ the
  function $\mb A$ defined by \eqref{eq:def:var-mass-embedding}
  satisfies the bounds
\begin{align*} \num\label{eq:var-mass-embedding-bound}
  &\|\mb A\|_{L^{p'}\L^{q'}(S_{m})}
    \lesssim \|\a\|_{L^{p'}\left(l^{r'}\right)}  &&p'\in\left(1,\,\infty\right]\;q'\in(r',\infty];
\end{align*}
furthermore the weak endpoint bounds 
\begin{align*}\num\label{eq:var-mass-embedding-bound-endpoint}
  &\|\mb A\|_{L^{p'}\L^{\infty}(S_{m})}\lesssim
    \|\a\|_{L^{p'}(l^{r'})} &&p'\in\left(1,\,\infty\right]\\
  &\|\mb A\|_{L^{1,\infty}\L^{q'}(S_{m})}\lesssim
    \|\a\|_{L^{1}(l^{r'})}&& q'\in(r',\infty]\\
  &\|\mb A\|_{L^{1,\infty}\L^{r',\infty}(S_{m})} 
    \lesssim \|\a\|_{L^{1}(l^{r'})}
\end{align*}
hold. All the above inequalities hold with constants independent of
the stopping sequence $\mf c$ appearing in
\eqref{eq:def:var-mass-embedding}.
\end{thm}
We refer
to Section \ref{sec:outer-measure-def}  for the description of the
outer measure structure on $\mb X$ and for the precise definition of
the iterated outer measure $L^{p}$ norms appearing on the left hand
sides.

\begin{cor}[Boundedness of the variational Carleson operator \cite{oberlin2012variation}]
  The operator \eqref{eq:var-carleson-op} defined pointwise for $f\in
  S(\R)$ extends to a bounded operator on $L^{p}(\R)$ for
  $r\in(2,\infty]$ and 
  $p\in(r',\infty)$.
\end{cor}

Given Theorem \ref{thm:max-var-mass-bounds} the above can be obtained
analogously as for the operator \eqref{eq:carleson-op}. For for $p$ and $r$ set
$\frac{1}{p'}=1-\frac{1}{p}$, $\frac{1}{r'}=1-\frac{1}{r}$, and choose $q$ and $q'$ so that
$\frac{1}{q}+\frac{1}{q'}=1$ and the bounds \eqref{eq:e-bounds} and
\eqref{eq:var-mass-embedding-bound} hold. Using the outer measure
Hölder inequality \eqref{eq:outer-duality} with the variational embedded function $\mb
A$ in lieu of $A$ and the wave-packet representation
\eqref{eq:var-duality-wavepacket} we obtain the required bound
\eqref{eq:var-duality-linear}.

Theorem \ref{thm:simple-mass-boundedness} follows from from Theorem
\ref{thm:max-var-mass-bounds} when $r=\infty$ by formally setting
\begin{align*}\num\label{eq:var-to-carleson-reduction}
&\a_{k}(z)=\left\{\begin{aligned}
    &a(z)&&\mbox{if }k=0\\
    &0&&\mbox{otherwise}
  \end{aligned}
\right.&&\mf c_{k}(z)=\left\{\begin{aligned}
    &-\infty&&\mbox{if } k<0\\
    &c(z)&& \mbox{if } k=0\\
    &+\infty && \mbox{if } k>0
  \end{aligned}
          \right.
\end{align*}

In particular the term $ \psi_{\eta,t}(y-z)\chi\big(t(\eta-c_{-})\big) $
appearing in \eqref{eq:A} are left truncated wave packets with respect
to the parameters $c_{-}$ and $c_{+}=+\infty$.


\subsection{Structure of the paper}\label{sec:structure-paper}

The rest of this paper is organized as follows.  In Section
\ref{sec:outer-measure-def} we define the outer measure structure on
$\mb X$.  We then recall properties of outer measure $L^{p}$ spaces and
generalize them to the iterated construction. In addition we
illustrate a limiting argument for maps to outer measure $L^{p}$ spaces
that allows to consider the bounds \eqref{eq:sm-bounds},
\eqref{eq:e-bounds}, and \eqref{eq:var-mass-embedding-bound} as
a-priori estimates. We also prove interpolation inequalities that
allow us to restrict the proof only to the the weak endpoints of the
above bounds.  Finally, we formulate the abstract outer Hölder
inequality and an outer Radon-Nikodym Lemma that imply inequality
\eqref{eq:outer-duality}.

In Section \ref{sec:var-carleson-embedding-reduction} we prove the
wave-packet domination bound \eqref{eq:var-duality-wavepacket}. In
particular it is shown that one can choose both the geometric
parameters of the outer measure space (see Section
\ref{sec:outer-measure-def}) and the parameters of the truncated
wave-packets in a compatible way i.e. so that both Thoerems
\ref{thm:energy-boundedness} and \ref{thm:max-var-mass-bounds} as well
as the conditions \eqref{eq:wave-packet-freq-support},
\eqref{eq:eq:wave-packet-geomtery-vanish} hold.  This is done by
providing a wave-packet representation for multipliers of the form
$\1_{[c_{-},c_{+})}$ with $c_{-}<c_{+}\in\R\cup\{+\infty\}$. For any
stopping sequence $\mf c$ this yields an embedded function
$\mb A_{\mf c}(y,\eta,t)$ so that
\begin{align*}\num\label{eq:var-duality-linearized}
  \int_{\R} \sum_{k\in\Z}\a_{k}(z)
  \int_{\c_{k}(z)}^{\c_{k+1}(z)} \FT f(\xi) e^{i\xi z}\dd \xi
  \;\dd z 
=
  \iiint_{\mb X} F(y,\eta,t)\mb A_{\mf c}(y,\eta,t) \dd y \dd \eta \dd t.
\end{align*}
The embedded function $\mb A_{\mf c}$ is pointwise dominated by
$\mb A$ and the map
$\a \mapsto \mb A_{\mf c}$ is shown to be
linear.  Furthermore the same procedure shows that the inequality in
\eqref{eq:duality-wave-packet} is actually an equality i.e.
\begin{align*}\num\label{eq:carleson-duality-linearized}
    \int_{\R} C_{c} f(z)\,a(z)\dd z = \iiint_{\R\times\R\times \R^{+}} F(y,\eta,t)A(y,\eta,t)\,
  \dd  y \dd \eta \dd t.
\end{align*}

In Section \ref{sec:aux-embedding} we introduce an auxiliary embedding
map for which we show iterated outer measure bounds. The crucial
result is given by the covering Lemma \ref{lem:superlevel-covering}
that allows one to control the measure of super-level sets of this
embedding map and by a projection Lemma \ref{lem:mass-projection} that
implies iterated bounds.

In Section \ref{sec:main-prop-proof} we actually prove Theorem
\ref{thm:max-var-mass-bounds} by showing the the auxiliary embedding
map of Section \ref{sec:aux-embedding} dominates the embedding
\eqref{eq:def:var-mass-embedding} in terms of sizes.

Finally, in Section \ref{sec:e-bounds} we show that bound
\eqref{eq:e-bounds} holds: this follows from an adaptation of the 
results of \cite{di2015modulation}. We also 
remark how in the case $p \in (2,r)$ a non-iterated version of outer
measure $L^{p}$ spaces is enough to obtain $L^{p}$ bounds for
\eqref{eq:var-carleson-op} and thus for \eqref{eq:carleson-op} with
$p\in(2,\infty)$.


\subsection{Notation}\label{sec:notation}

We quickly recall some useful notation.

We say that $A(x)\lesssim B(x)$ if there exists a constant $C>0$ such
that $A(x)\leq C B(y)$ for all $x,y$ in the domains of $A$ and $B$
respectively. Unless otherwise specified the constant $C>0$ is
absolute. We may emphasize the dependence on a specific parameter $p$ by
writing $A(x)\lesssim_{p} B(y)$. We write $A(x)\approx B(y)$ if
$A(x)\lesssim B(y)$ and $A(x)\gtrsim B(y)$.

We denote open and close Euclidean balls of $\R $ as 
\begin{align*}
&  B_{r}(x) :=(x-r,x+r) && B_{r}:= (-r,+r) &&  \closure{B_{r}(x)} :=[x-r,x+r] && \closure{B_{r}}:= [-r,+r].
\end{align*}
We indicate by $\1_{\Theta}$ the characteristic function of the set
$\Theta$ i.e.
\begin{align*}
  \1_{\Theta}(x) :=\left\{\begin{aligned}
      &      1&&\text{if } x\in\Theta\\
      &      0&&\text{if } x\notin \Theta
    \end{aligned}
  \right.
\end{align*}
For an arbitrary large
$N>0$ we introduce the smooth bump function
\begin{align*}\num\label{eq:def:bump}
&  W(z) := \left(1+|z|^{2}\right)^{-N/2} && W_{t}(z) := t^{-1}W\left(\frac{z}{t}\right).
\end{align*}
We define
\begin{align*}
\fint_{B_{r}(x)}f(z)\dd z := \frac{1}{2r} \int_{B_{r}(x)} f(z)\dd z.
\end{align*}
The operators  $M$ and $M_{p}$ are the Hardy-Littlewood maximal function i.e.
    \begin{align*}\num\label{eq:HL-maximal-function}
      &Mf(z) :=    \sup_{t\in\R^{+}}\fint_{B_{t}(z)}|f(z')|\dd z'\\
      &M_{p}f(z) :=    \sup_{t\in\R^{+}}\Big(\fint_{B_{t}(z)}|f(z')|^{p}\dd z'\Big)^{1/p}.
    \end{align*}    
Given a function $\phi\in S(R)$ we obtain its  frequency translates
and dilates by setting
\begin{align*}
  \phi_{\eta,t}(z):=t^{-1}e^{i\eta z} \phi\left(\frac{z}{t}\right).
\end{align*}
The stopping sequence $\c$ will denote  a Borel measurable function defined on
$\R$ with values in increasing sequences in $\R\cup {+\infty}$ i.e.
\begin{align*}
&  z\mapsto \c(z)=(\c_{k}(z))_{k\in\Z}& &-\infty<\dots\leq\c_{k-1}(z)\leq\c_{k}(z)\leq \c_{k+1}(z)\leq \dots\leq +\infty.
\end{align*}
Similarly $\a$ will denote a Borel Bochner-measurable function on $\R$ with
values in $l^{r'}$ i.e.
\begin{align*}
&  z\mapsto \a(z)=(\a_{k}(z))_{k\in\Z}\in l^{r'}.
\end{align*}

We use the notation $L^{p}(S)$ and $L^{p}\L^{q}(S)$ to denote (iterated)
outer measure $L^{p}$ spaces. The (outer-) measure of the space is
omitted from the notation.  We distinguish the above from $L^{p}$ that
are classical Lebesgue spaces. In the case of of $L^{p}$ spaces on
$\R$ the measure is the Lebesgue measure; when necessary we may
emphasize the measure $\mc L$ on the space by writing $L^{p}_{(\dd\mc L)}$.


\section{Outer measures on the time-frequency
  space}\label{sec:outer-measure-def}

We begin the description of the outer measure on the time-frequency
space $\mb X$ by introducing a family of distinguished generating
sets.  The \emph{tent} $T(x,\xi,s)\subset \mb X$ indexed by the top
point $(x,\xi,s)\in\mb X$ is the set
\begin{align*}\num\label{eq:def:tents}
  &T(x,\xi,s) := T^{(i)}(x,\xi,s) \cup T^{(e)}(x,\xi,s)
    \\
  & T^{(i)}(x,\xi,s) := \left\{ (y,\eta,t) \sthat |y-x|<s,\,t(\eta-\xi) \in
    \Theta^{(i)},\; t<s\right\} \\
  & T^{(e)}(x,\xi,s) := \left\{ (y,\eta,t) \sthat |y-x|<s,\,
    t(\eta-\xi)\in \Theta^{(e)},\; t<s\right\}
\end{align*}
where
\begin{align}\label{eq:geometric-intervals}
    \Theta=(\alpha^{-},\alpha^{+})&&\Theta^{(i)}=(\beta^{-},\beta^{+})
    &&\Theta^{(e)}=\Theta\setminus\Theta^{(i)}
  \end{align}
  are geometric intervals such that
  $0\in\Theta^{(i)}\subseteq\Theta$
  i.e. $\alpha^{-}\leq\beta^{-}<0<\alpha^{+}\leq \beta^{+}$.  We refer
  to $T^{(i)}$ and $T^{(e)}$ as the interior and exterior parts of the
  tent $T$. To define the iterated outer measure structure we
  introduce \emph{strips} $D(x,s)\subset \mb X$ as
  \begin{align}\label{eq:def:stripes}
    &D(x,s) := \left\{ (y,\eta,t)\;\sthat\;|y-x|<s,\, t<s\right\}.
  \end{align}
  We indicate the family of all tents by $\mb T$ and the family of all
  strips by $\mb D$.

  \begin{figure}[h]
    \centering
    \begin{tikzpicture}[scale=3]

      \coordinate (o) at (0,0); \coordinate (xmax) at (1.5,0);
      \coordinate (xmin) at (-0.1,0); \coordinate (ymax) at (0,1);
      \coordinate (ymin) at (0,0);

      \draw[->] (xmin) -> (xmax); \draw[->] (ymin) -> (ymax);

      \coordinate (top) at (0.7,0.9);

      \node [left] at (ymax) {$t$}; \node [below] at (xmax) {$y$};

      \node [above] at (top) {$(x,s)$}; \draw [fill] (top) circle
      (0.02);

      \draw let \p1 = (top) in node [anchor = north ] at
      ($(\x1,0)+0.5*(\y1,0)$) {$x+s$}; \draw let \p1 = (top) in [fill]
      ($(\x1,0)+0.5*(\y1,0)$) circle (0.02);

      \draw let \p1 = (top) in node [anchor = north ] at
      ($(\x1,0)-0.5*(\y1,0)$) {$x-s$}; \draw let \p1 = (top) in [fill]
      ($(\x1,0)-0.5*(\y1,0)$) circle (0.02);

      \draw let \p1 = (top) in 
      (top) -- ($(top)+0.5*(\y1,0)$) -- ($(\x1,0)+0.5*(\y1,0)$) --
      ($(\x1,0)-0.5*(\y1,0)$) --($(top)-0.5*(\y1,0)$)--cycle;

    \end{tikzpicture}
    \begin{tikzpicture}[scale=3]

      \coordinate (o) at (0,0); \coordinate (xmax) at (1.5,0);
      \coordinate (xmin) at (-0.1,0); \coordinate (ymax) at (0,1);
      \coordinate (ymin) at (0,0);

      \draw[->] (xmin) -> (xmax); \draw[->] (ymin) -> (ymax);

      \coordinate (top) at (0.7,0.9);

      \node [anchor= north east] at (ymax) {$t$}; \node [below] at
      (xmax) {$\eta$};

      \node [above] at (0.7,0.95) {$(\xi,s)$}; \draw [fill] (top)
      circle (0.02);

      \draw [domain= 1:1.5, variable =\x] plot
      ({\x},{0.9*0.3/(\x-0.7)}); \draw [domain= -0.1:0.5, variable
      =\x] plot ({\x},{-0.9*0.2/(\x-0.7)});

      \draw [dotted, domain= 0.8:1.5, variable =\x] plot
      ({\x},{0.9*0.1/(\x-0.7)}); \draw [dotted, domain= -0.1:0.6,
      variable =\x] plot ({\x},{-0.9*0.1/(\x-0.7)});

      \draw (0.5,0.9) -- (1,0.9);

      \draw node [anchor =south west ] at (1,0.9) {\tiny
        $\xi\!\!+\!\alpha^{\!+}\!s^{-1}$}; \draw [fill] (1,0.9) circle
      (0.02);

      \draw node [anchor =south] at (0.5,0.9) {\tiny
        $\xi\!\!+\!\alpha^{\!-}\!s^{-1}\hspace{2.8em}$}; \draw [fill]
      (0.5,0.9) circle (0.02);

      \draw node at (0.7,0.2) {$T^{(i)}$};
      \draw node at (1.1,0.50) {$T^{(e)}$};
      \draw node at (0.4,0.4) {$T^{(e)}$};

    \end{tikzpicture}
    
    \caption{The tent $T(x,\xi,s)$.}
    \label{fig:tent}
  \end{figure}
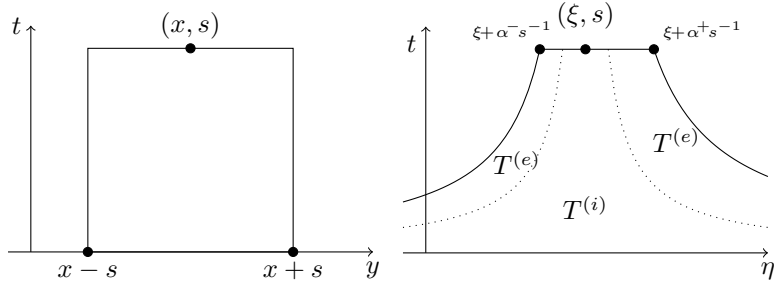
  
  The specific values of the geometric intervals $\Theta$,
  $\Theta^{(i)}$, and $\Theta^{(e)}$ in \eqref{eq:def:tents} are often
  inessential.  However, the freedom of choosing appropriate
  parameters was shown to be important in \cite{do2015lp}.  Theorem
  \ref{thm:max-var-mass-bounds} holds as long as
\begin{align*}\num\label{eq:truncate-geometry}
  &\closure{B_{b}}\subset \Theta^{(i)}\subset \closure{B_{d'}}
  &&\closure{B_{\epsilon}(d)}\cup \closure{B_{\epsilon}(-d)}\subset\Theta^{(e)}
\end{align*}
with $b$, $d$, $d'$, and $\epsilon$ appearing in
\eqref{eq:wave-packet-freq-support},
\eqref{eq:eq:wave-packet-geomtery-vanish}, and
\eqref{eq:eq:wave-packet-geomtery-constant}. As a matter of fact, if
one were to consider only left truncated wave-packets in
\eqref{eq:def:var-mass-embedding} then Theorem
\ref{thm:max-var-mass-bounds} would hold as long as
\begin{align*}\num\label{eq:truncate-geometryl}
  &B_{b}\subset \Theta^{(i)}&& -d'<\beta^{-} &&\closure{ B_{\epsilon}(d)}
                                             \subset \Theta^{(e)}\cap \R^{+}=[\beta^{+},\alpha^{+}).
\end{align*}
Theorem \ref{thm:simple-mass-boundedness} holds as long as satisfies
\begin{align*}\num\label{eq:chi-geometry}
  &B_{b}\subset \Theta^{(i)} &&\spt \chi\subset \Theta^{(e)}\cap \R^{+}= [\beta^{+},\alpha^{+}).
\end{align*}
Theorem \ref{thm:energy-boundedness} holds as long as
$  B_{b}\subset\Theta^{(i)}$.  From now on we will allow all our
implicit constants to depend on $\Theta$ and $\Theta^{(i)}$.

We now define the outer measures $\mu$ and $\nu$ by introducing the
pre-measures $\overline{\mu}$, and $\overline{\nu}$ on the generating
sets 
\begin{align*}\num\label{eq:def:premeasures}
&  \overline{\mu}(T(x,\xi,s)):=s&&  \overline{\nu}(D(x,s)):=s.
\end{align*}
The outer measure of an arbitrary
subset $E\subset \mb X$ are obtained via a covering procedure using
\emph{countable} unions of generating sets i.e.
\begin{align*}\num\label{eq:def:outer-measures}
  &\nu(E) := \inf\Big\{ \sum_{n\in\N} \overline{\nu}(D_{n}) \sthat 
    E\subset \bigcup_{\substack{D_{n}\in\mb D\\n\in\N}} D_{n},\Big\}
\end{align*}
and similarly for $\mu$ using $\overline{\mu}$ and the family $\mb
T$. We say that $\nu$ and $\mu$ are generated by the pre-measures
$(\overline{\nu},\mb D)$ and $(\overline{\mu},\mb T)$ respectively. We
call an \emph{outer measure space} a pair $(\mb X,\mu)$ of a
separable complete measure space $\mb X$ and an outer measure
$\mu:2^{\mb X}\to \R^{+}\cup\{+\infty\}$. We
will henceforth suppose that $\mu$ is generated by pre-measures
$(\overline{\mu},\mb T)$ where $\mb T$ is a collection of subsets
$T\subset \mb X$ that we assume to be Borel measurable.

The final ingredient we need for introducing outer measure $L^{p}$
spaces is a notion of how large a function on $\mb X$ is.  We call a
\emph{size} any quasi-norm $\|\cdot\|_{S}$ on Borel functions on
$\mb X$ i.e. a positive functional that satisfies the following
properties.
{\bf Monotonicity:} for any Borel function $G_{1}$ and
$G_{2}$
\begin{align*}\num\label{eq:size-monotonicity}
  |G_{1}|\leq|G_{2}|\implies \|G_{1}\|_{S}\lesssim
    \|G_{2}\|_{S}.
\end{align*}
{\bf Positive homogeneity:} for all Borel functions
    $G$  
    \begin{align*}\num\label{eq:size-homogeneity}
          \|\lambda G \|_{S}=|\lambda| \|G\|_{S}  \qquad \forall \lambda\in \C.
    \end{align*}
    {\bf Quasi-triangle inequality:} for any sequence of Borel functions
    $G_{k}$ and for some quasi-triangle constant
    $c_{s}\geq 1$
    \begin{align*}\num\label{eq:size-triangle}
      \|\sum_{k=0}^{\infty}G_{k}\|_{S}\leq \sum_{k=0}^{\infty}
    c_{s}^{k+1}\|G_{k}\|_{S} 
    \end{align*}
    
    We define the $S,\mu$ -  super-level outer measure as
\begin{align}
  \mu\left(
  \|G\|_{S}>\lambda
  \right) := \inf \left\{\mu(E_{\lambda})\,
  \sthat \, \|G\,\1_{\mb X \setminus E_{\lambda}}
  \|_{S}\leq\lambda\right\}
\end{align}
where the lower bound is taken over Borel subset $E_{\lambda}$ of
$\mb X$. The outer-$L^{p}$ quasi-norms for $p\in(0,\infty]$ are give by
\begin{align*}
  \left\| G \right\|_{L^{p}\left(S\right)}^{p} :=
  \int_{\lambda\in\R^{+}} p \lambda^{p} \mu \left(\|G\|_{S} >
    \lambda\right)\frac{\dd \lambda}{\lambda};
\end{align*}
weak outer $L^{p}$ quasi-norms are similarly given by
\begin{align*}
  \left\| G \right\|_{L^{p,\infty}(S)}^{p} :=  \sup_{\lambda\in \R^{+}} p
  \lambda^{p}\mu \left(\|G\|_{S} > \lambda\right).
\end{align*}
The outer $L^{p}$ spaces are subspaces of Borel functions on $\mb X$
for which the above norms are finite.  The expressions
defining outer $L^{p}$ quasi-norms are based on the super-level set
representation of the Lebesgue integral, however the expression
$\mu \left(\|G\|_{S} > \lambda\right)$ that appears in lieu of the classical
$\mu\left(\{x \,\sthat\, |g(x)|>\lambda\}\right)$  cannot always be
interpreted as a measure of a specific set.  Generally speaking,
$L^{p}$ spaces for $p\in(0,\infty)$ are interpolation spaces between
the size quasi-norm and the outer measure of the support of a function.

Using a slight abuse of notation we say that a size $\|\cdot\|_{S}$ is
generated by $(\|\cdot\|_{S(T)},\mb T)$ where $\|\cdot\|_{S(T)}$ are
sizes indexed by generating sets $T\in\mb T$ and in particular
\begin{align*}\num\label{eq:size-generation}
  \|G\|_{S}:=\sup_{T\in\mb T} \|G\|_{S(T)}.
\end{align*}

The construction of iterated outer $L^{p}$ spaces is based on using
localized versions of outer $L^{q}$ quasi-norms as sizes themselves.
Notice that outer $L^{q}$ norms are quasi-norms since they too satisfy
the quasi-triangle inequality.  Given a size $S$ and a generating
pre-measure $(\overline{\nu},\mb D)$, outer $\L^{q}(S)$ sizes are
generated by $(\L^{q}(S)(D),\mb D)$ where 
\begin{align*}\num\label{eq:def:L-size}
&    \|G\|_{\L^{q}(S)(D)}:= \frac{\left\|F\1_{D}\right\|_{L^{q}(S)}}{\nu \left( D \right)^{1/{q}}}
\end{align*}
so $\|G\|_{\L^{q}(S)} := \sup_{D\in\mb D}\|G\|_{\L^{q}(S)(D)}$.
Consequently we construct iterated outer $L^{p}$ spaces as
\begin{align*}\num\label{eq:def:LpLqS}
  \|G\|_{L^{p}\L^{q}(S)}^{p}:=\int_{\tau\in\R^{+}}p\,
  \tau^{p}\,\nu\left(\|G\|_{\L^{q}(S)}>\tau\right)\frac{\dd \tau}{\tau}.
\end{align*}

To deal with embedded functions $F$ and $\mb A$ from \eqref{eq:F} and
\eqref{eq:def:var-mass-embedding} we introduce the respective sizes
$\|\cdot\|_{S_{e}}$ and $\|\cdot\|_{S_{m}}$ that are generated by
$(S_{e}(T),\mb T)$ and $(S_{m}(T),\mb T)$ respectively. The two
families of ``local'' sizes $S_{e}(T)$ and $S_{m}(T)$ are given by
\begin{align*}    \num\label{eq:def:size-e}
  &\begin{aligned}
    \|F\|_{S_{e} (T)} := &\;  \frac{ \left\|F\1_{T^{(e)}}\right\|_{L^{2}}}{\mu \left( T \right)^{1/2}}
    + \left\|F\1_{T}\right\|_{L^{\infty}}\\
    =&\; \|F\|_{S^{2}(T^{(e)})} +\|F\|_{S^{\infty}(T)}
  \end{aligned}\\
  \num\label{eq:def:size-m}
  &\begin{aligned}
    \|\mb A\|_{S_{m}(T)} := &\; \frac{      \left\|\mb A\1_{T}\right\|_{L^{2}}}{\mu \left( T \right)^{1/2}}
+ \frac{\left\|\mb A\1_{T^{(i)}}\right\|_{L^{1}} }{\mu \left( T \right)} 
    \\=&\;\|\mb A\|_{S^{2}(T)} +\|\mb A\|_{S^{1}(T^{(i)})}.
\end{aligned}
\end{align*}
Here  $L^{2}$, $L^{\infty}$, and $L^{1}$ refer to classical Lebesgue
$L^{p}$ norms on $\mb X$ with respect to the Borel measure
$\dd y\dd \eta\dd t$.  The local sizes $\|\cdot\|_{S_{e}(T)}$ coincide with the ones introduced on
the upper 3-space in \cite{do2015lp} while $\|\cdot\|_{S_{m}(T)}$ are
dual to the former in an appropriate sense. 

We conclude the construction of outer measure $L^{p}$ spaces with a
useful remark about the specific geometric properties of coverings
with tents $\mb T$. For any tent $T(x,\xi,s)$ we define its 
$R$-enlargement with $R>1$ as
\begin{align*}\num\label{eq:tree-enlargement}
  &R\,T(x,\xi,s) := \bigcup_{{|\xi'-\xi|<R s^{-1}}} T(x,\xi',Rs).
\end{align*}
Notice that $\mu(R\,T)\lesssim R^{3} \mu(T)$ with a constant that
depends on the geometric intervals \eqref{eq:geometric-intervals} but
not on $R$. As a matter of fact the set $R\, T$ can be covered by a
finite collection of tents $T(x,\xi_{i},Rs)$ by choosing $\xi_{i}$
such that
  \begin{align*}
  \bigcup_{i} \left\{\eta \sthat Rs(\eta-\xi_{i})\in \Theta\right\}
  \supset \bigcup_{|\xi'-\xi|<Rs^{-1}} \left\{\eta \sthat Rs(\eta-\xi')\in \Theta\right\}
  \end{align*}
  The number of points $\xi_{i}$ needed to do this is bounded up to a
  constant factor by $R^{2}$ and thus
  $\mu\left(R\,T\right)\lesssim R^{3} \mu(T)$.


\subsection{Properties of outer measure $L^{p}$
  spaces.}\label{sec:outer-measure-properties}

We recall some important properties of outer measure $L^{p}$ spaces and
elaborate on how they carry over to iterated outer-measure
spaces. Generally $\mb X$ may be any locally compact complete metric
space; in our case $\mb X=\R\times\R\times \R^{+}$ with
\begin{align*}
\dist\big((y,\eta,t);(y',\eta',t')\big)=
t^{-1}|y-y'|+t|\eta-\eta'|+ |\log\frac{t}{t'}|.
\end{align*}

\subsubsection{Dominated convergence}\label{sec:outer-dominated-convergence}
While outer measure $L^{p}$ spaces fall into the class of quasi-Banach
spaces, we record only some functional properties that are useful for 
our applications.

Recall that the quasi-triangle inequality for sizes
\eqref{eq:size-triangle} holds for both
finite and infinite sums.  Given an outer measure space $(\mb X,\mu)$
and a size $\|\cdot\|_{S}$, the outer measure $L^{p}$ quasi-norms also
satisfy the quasi-triangle inequality:
\begin{align*}\num\label{eq:triangle-series}
    \Big\|\sum_{k=0}^{\infty}G_{k}\Big\|_{L^{p}(S)}\lesssim_{c_{s}',c_{s},p}
    \sum_{k=0}^{\infty} c_{s}'^{k+1}\|G_{k}\|_{L^{p}(S)}
  \end{align*}
  for any $c_{s}'> c_{s}$ where $c_{s}$ is the quasi-triangle constant
  of the size $S$. As a matter of fact, for any $\lambda>0$ and for
  every $k$ choose $E_{\lambda,k}$ such that
  \begin{align*}
    &    \|G_{k}\1_{\mb X\setminus E_{\lambda,k}}\|_{S}\leq \lambda &
    &
      \|G_{k}\|_{L^{p}(S)}^{p} \lesssim
      \int_{\R^{+}} p \lambda ^{p} \mu(E_{\lambda,k})\frac{\dd \lambda}{\lambda}
  \end{align*}
  and set
  $ E_{\lambda}=\bigcup_{k=0}^{\infty}E_{\lambda c_{s}'^{-k-1},k}$ so
  that using the quasi-triangle inequality for $\|\cdot\|_{S}$ one has 
    \begin{align*}
      &\mu(E_{\lambda})\leq\sum_{k=0}^{\infty}\mu\left(E_{\lambda
        c_{s}'^{-k-1},k}\right) &
      &
        \Big\|\sum_{k=0}^{\infty}G_{k} \1_{E_{\lambda}^{c}}\Big\|_{S}\leq  \lambda\frac{c_{s}}{(c_{s}'-c_{s})}.
    \end{align*}
    Thus
    \begin{align*}
      \Big\|\sum_{k=0}^{\infty}G_{k}\Big\|_{L^{p}(S)}^{p}&
                                                                  \leq 
      p\left(\frac{c_{s}}{c_{s}'-c_{s}}\right)^{p}\int_{\R^{+}}\lambda^{p}\mu(E_{\lambda})\frac{\dd
      \lambda}{\lambda}\\&\leq 
      \left(\frac{c_{s}}{c_{s}'-c_{s}}\right)^{p}\sum_{k=0}^{\infty} 
                           c_{s}'^{p(k+1)} \|G_{k}\|_{L^{p}(S)}^{p}.
    \end{align*}
    If $p\geq 1$ then this concludes the proof. Otherwise for any
    $\epsilon>0$ one has
    \begin{align*}
      \left\|\sum_{k=0}^{\infty}G_{k}\right\|_{L^{p}(S)} \lesssim_{\epsilon,p,c,c_{s}} \sum_{k=0}^{\infty}(1+\epsilon)^{k+1}c_{s}'^{k+1}\|G_{k}\|_{L^{p}(S)}
    \end{align*}
    but since $c_{s}'>c_{s}$ was arbitrary this also allows us to
    conclude.
    
    This fact is crucial to be able to use localized
    outer $\L^{p}$ quasi-norms as sizes themselves. Furthermore we
    deduce the following  domination property.
  \begin{cor}\label{cor:triangle-domination}
    Suppose that $G$ is a Borel function on $\mb X$ and  $|G|\leq \limsup_{n\to \infty}|G_{n}|$ pointwise on
    $\mb X$ for some sequence of Borel functions $G_{n}$ that satisfy
  \begin{align*}
&    \|G_{n+1}-G_{n}\|_{L^{p}(S)}\leq C\,
    c_{s}'^{-n}\|G_{0}\|_{L^{p}(S)}&&\text{for some }  c_{s}'> c_{s}.
                                 \shortintertext{Then}
  &                               \|G\|_{L^{p}(S)}\lesssim_{C,p,c_{s}',c_{s}}
  \|G_{0}\|_{L^{p}(S)}.\end{align*}  
\end{cor}
This follows from \eqref{eq:triangle-series} and from the
monotonicity properties of sizes and thus of outer $L^{p}$ quasi-norms.

Using this property we will restrict ourselves to proving bounds
\eqref{eq:sm-bounds}, \eqref{eq:e-bounds}, and
\eqref{eq:var-mass-embedding-bound} for a dense class of functions. In
particular we will always consider the functions in play to be smooth
and rapidly decaying.  For example, given a function
$\a\in L^{p}(l^{r'})$ one may always choose a sequence of
approximating functions
$ \a^{(n)}\subset C_{c}^{\infty}(l^{r'})$ such that
\begin{align*}
  &\|\a^{(0)}\|_{L^{p}(l^{r'})}\lesssim\|\a\|_{L^{p}(l^{r'})}\\
  &\|\a^{(n+1)}-\a^{(n)}\|_{L^{p}(l^{r'})}\lesssim 2^{-Nn}\|\a\|_{L^{p}(l^{r'})}
\end{align*}
for an arbitrary $N>1$. Considering the sequence embedded functions
$\mb A_{n}$ associated to $\a^{(n)}$ via
\eqref{eq:def:var-mass-embedding}, the pointwise relation
$\mb A=\lim_{n}\mb A_{n}$ clearly holds.  Corollary
\ref{cor:triangle-domination} applied to $\mb A_{n}$ allows us to
conclude that if bounds of Theorem \ref{thm:max-var-mass-bounds} hold for the
functions $\a^{(n)}$ they also hold for $\a$. Thus we can restrict to
  proving the bounds as a priori estimates i.e. we can restrict to
  showing that they hold for a dense class of functions $\a$. The same
  can be done for the energy embedding bounds of
  Theorem \ref{thm:energy-boundedness}.

\subsubsection{Hölder and Radon-Nikodym inequalities}

We now illustrate the abstract outer measure results from which
inequality \eqref{eq:outer-duality} follows. The first two statements
relate to general outer measure spaces and are similar to what was
obtained in \cite{do2015lp}.

\begin{lem}[Radon-Nikodym domination]\label{lem:radon-nikodym}
  Consider $(\mb X,\mu)$ an outer measure space with $\mu$, generated
  by $(\overline{\mu},\mb T)$ as in \eqref{eq:def:outer-measures},
  endowed with a size $\|\cdot\|_{S}$ generated by
  $(\|\cdot\|_{S(T)},\mb T)$.  Suppose that the generating family
  $\mb T$ consists of Borel sets and satisfies the \emph{covering
    condition} i.e.  $\mb X =\bigcup_{i\in\N}T_{i}$ for some countable
  sub-collection $T_{i}\in\mb T$.
  
  If $\mc L$ is a positive Borel measure on $\mb X$ such that 
  \begin{align*}\num\label{eq:rn-diff}
    &    \int_{T}|G(P)|\dd \mc L(P)\leq C \|G\|_{S(T)}\,\overline{\mu}(T)
    & &\forall        T\in\mb T 
  \end{align*}
  and for any Borel function $G$ and 
  \begin{align*}\num\label{eq:rn-measure-absolutely-continuous}
    &    \mu(E)=0 \implies \mc L(E)=0 && \forall E\subset \mb X \text{ Borel}
  \end{align*}
  then for any Borel function $G$ the bound
  \begin{equation*}\num
    \Big|\int_{\mb X} G(P) \dd \mc L(P)\Big| \lesssim \big\| G
    \big\|_{L^{1}(S)}
\end{equation*}
holds.
\end{lem}
The proof of this Lemma is similar to the one in \cite{do2015lp}.
\begin{proof}
  Suppose $\|G\|_{L^{1}(S)}<\infty$, otherwise there is
  nothing to prove.  For each $k\in\Z$ let $E_{2^{k}}'$ be a Borel set
  such that
  \begin{align*}
    &\|G \1_{\mb X \setminus E_{2^{k}}'}\|\leq 2^{k}& \mu(E_{2^{k}}') \leq 2
                                                  \mu\left(\|G\|_{S}>2^{k}\right).
  \end{align*}
  so $\|G\|_{L^{1}(S)}\lesssim \sum_{k=-\infty}^{+\infty} 2^{k}\mu(E_{2^{k}}')$.
  Set
  \begin{align*}
    &    E_{2^{k}}:=\bigcup_{l=k}^{+\infty} E_{2^{k}}' &&\Delta E_{2^{k}}:=
                                                    E_{2^{k-1}}\setminus
                                                    E_{2^{k}} && E_{0}=\bigcup_{k=-\infty}^{+\infty}E_{2^{k}}&&E_{\infty}=\bigcap_{k=-\infty}^{+\infty}E_{2^{k}}.
  \end{align*}
  We have
  \begin{align*}
    \Big|\int_{\mb X} G(P)\dd \mc L(P)\Big|\leq \int\mcl_{\mb X\setminus E_{0}}|G(P)|\dd \mc L(P) + \sum_{k=-\infty}^{+\infty} \int\mcl_{\quad\Delta
    E_{2^{k}}} |G(P)|\dd \mc L(P)  + \int\mcl_{E_{\infty}} |G(P)|\dd \mc L(P).
  \end{align*}
  where
  \begin{align*}
    &\|G \1_{\Delta E_{2^{k}}}\|_{S} \leq 2^{k}&
    &    \|G\|_{L^{1}(S)} \lesssim \sum_{k=-\infty}^{+\infty} 2^{k} \mu(\Delta E_{2^{k}}).
  \end{align*}
  For every $k$ there exists a countable covering
  $\bigcup_{l\in\N}T_{k,l}\supset \Delta E_{2^{k}}$ such that
  \begin{align*}
  \sum_{l\in\N}\overline{\mu}(T_{k,l})\leq 2 \mu(\Delta
  E_{2^{k}}).
  \end{align*}
  For each $k\in\Z$ apply \eqref{eq:rn-diff} to obtain
  \begin{align*}
    \int_{\Delta E_{2^{k}}}|G(P)|\dd \mc L(P) &\leq \sum_{l\in\N} \int_{T_{k,l}}
    |G(P)|\1_{\Delta E_{2^{k}}}(P)\dd \mc L(P)\\
 &\leq \|G\1_{\Delta E_{2^{k}}}\|_{S}\sum_{l\in\N}\overline{\mu}(T_{k,l}) \leq 
   2^{k+1}\mu(\Delta E_{2^{k}})
  \end{align*}
  Thus
  \begin{align*}
    \sum_{k=-\infty}^{+\infty} \int_{\Delta E_{2^{k}}} |G(P)|\dd \mc L(P)\lesssim \|G\|_{L^{1}(\mb X,\mu,S)}.
  \end{align*}
  The term $\int_{{\mb X\setminus E_{0}}}|G(P)|\dd \mc L(P) $
  vanishes because we may represent $\mb X=\bigcup_{i\in\N}
  T_{i}$. Using \eqref{eq:rn-diff} and the monotonicity of sizes we
  have
  \begin{align*}
    \int_{\mb X\setminus E_{0}}|G(P)|\dd \mc L
   (P)&\leq\sum_{i\in\N}
    \int_{T_{i}}|G(P)|\1_{\mb X\setminus E_{0}}(P) \dd \mc
  L(P)
    \\ &\lesssim\sum_{i\in\N}
    \|G\1_{\mb
    X\setminus E_{0}}\|_{S(T_{i})} \,\overline{\mu}(T_{i})=0.
  \end{align*}
  The term $\int_{{ E_{+\infty}}} |G(P)|\dd \mc L(P)$ also vanishes
  since
  \begin{align*}
  \mu(E_{2^{k}})\le \sum_{l=k}^{\infty}\mu(E_{2^{k}}')\lesssim
  2^{-k}\|G\|_{L^{1}(S)}
  \end{align*}
and thus $\mu(E_{+\infty})=0$ and  $\mc L(E_{+\infty})=0$ by
  \eqref{eq:rn-measure-absolutely-continuous}. This concludes the
  proof.
\end{proof}

The proof of the following outer measure Hölder inequality can be
found in \cite{do2015lp}.

\begin{prop}[Outer Hölder inequality]\label{prop:outer-holder}
  Let $(\mb X,\mu)$ be an outer measure space endowed with three sizes
  $\|\cdot\|_{S}$, $\|\cdot\|_{S'}$, and $\|\cdot\|_{S''}$ such that
  for any Borel functions $F$ and $A$ on $\mb X$ the product
  estimate for sizes
  \begin{align*}\num\label{eq:holder-size-condition}
    \|FA\|_{S}\lesssim \|F\|_{S'}\|A\|_{S''}
  \end{align*}
  holds.  Then for any Borel functions $F$ and $A$ on $\mb X$ the
  following outer Hölder inequality holds:
\begin{align*}\num\label{eq:outer-holder}
\|F  A\|_{L^{p}(S)}\leq 2 \|F\|_{L^{p'}(S')}\| A\|_{L^{p''}(S'')}
\end{align*}
for any triple  $p,p',p''\in(0,\,\infty]$ of exponents such that $\frac{1}{p'}+\frac{1}{p''}=\frac{1}{p}$,
\end{prop}

The above two statement can be easily extended to iterated outer
measure spaces. Suppose from now on that $\mb X$ is endowed with two
outer measures $\nu$ and $\mu$, the former generated by a pre-measure
$(\overline{\nu},\mb D)$ as described in
\eqref{eq:def:outer-measures}. Given a size $\|\cdot\|_{S}$ we
introduce local $\L^{q}(S)$ sizes as described by
\eqref{eq:def:L-size} and the corresponding iterated outer $L^{p}\L^{q}(S)$
quasi-norms as described in \eqref{eq:def:LpLqS}.

\begin{cor}[Outer Hölder inequality for iterated outer measure spaces]\label{cor:iter-holder}
  Let $(\mb X,\mu)$ be an outer measure space endowed with three sizes
  $\|\cdot\|_{S}$, $\|\cdot\|_{S'}$, and $\|\cdot\|_{S''}$ satisfying the
  assumptions of Proposition \ref{prop:outer-holder}. Then given any
  two triples pairs of  exponents $p,p',p''\in(0,\,\infty]$ and
  $q,q',q''\in(0,\infty]$ such that $\frac{1}{p'}+\frac{1}{p''}=\frac{1}{p}$ and
  $\frac{1}{q'}+\frac{1}{q''}=\frac{1}{q}$ the iterated Hölder bounds
  \begin{align*}\num\label{eq:iter-outer-holder}
    \|F A\|_{L^{p}\L^{q}(S)}\lesssim     \|F \|_{L^{p'}\L^{q'}(S')}    \|A\|_{L^{p''}\L^{q''}(S'')}
  \end{align*}
  hold for any Borel functions $F$ and $A$ on $\mb X$.
\end{cor}
As a matter of fact the inequality 
\begin{align*}\num\label{eq:iterated-holder-size}
  \|FA\|_{\L^{q}(S)}\lesssim   \|F\|_{\L^{q'}(S')}
  \|A\|_{\L^{q''}(S'')}
\end{align*}
holds for localized $\L^{q}(S)$ sizes satisfy the inequality by
 Proposition \ref{prop:outer-holder} applied to the defining expression
\eqref{eq:def:L-size}. Thus the local $\L^{q}(S)$ sizes themselves satisfy the conditions
of Hölder inequality and the statement of the above Corollary follows.

The Radon-Nikodym Lemma \ref{lem:radon-nikodym} can also be
generalized to iterated outer measure $L^{p}$ spaces.

\begin{cor}[Iterated Radon-Nikodym domination]\label{cor:iter-radon-nikodym}
  Consider $(\mb X,\mu)$ an outer measure space with a size
  $\|\cdot\|_{S}$ and a Borel measure $\mc L$ that satisfy the
  conditions of Lemma \ref{lem:radon-nikodym} and let $\nu$ be a
  measure generated by $(\overline{\nu},\mb D)$. Suppose that $\mb D$
  also satisfies the covering condition of Lemma
  \ref{lem:radon-nikodym}.  Then the iterated Radon-Nikodym domination
  \begin{align*}
    \Big|\int_{\mb X} G(P) \dd \mc L(P)\Big| \lesssim     \|G(P)\|_{L^{1}\L^{1}(S)}
  \end{align*}
  holds.
\end{cor}

As a matter of fact  for any Borel function $G$  the inequality 
\begin{align*}
  \int_{D} |G(P)|\mc L(P) \lesssim \|G\|_{\L^{1}(S)(D)} \overline{\nu}(D)
\end{align*}
follows from \eqref{eq:def:L-size} and Lemma
\ref{lem:radon-nikodym}. Thus the outer measure space $(\mb X,\nu)$
and the family of local sizes
$\|\cdot\|_{\L^{1}(S)(D)}$ satisfy the conditions of
Lemma \ref{lem:radon-nikodym} and the statement of the Corollary follows.

Using the above properties one can deduce inequality
\eqref{eq:outer-duality}: introduce the size
\begin{align*}
  \|G\|_{S^{1}}:=\sup_{T\in \mb T} \|G\|_{S^{1}(T)} =\sup_{T\in\mb T}
  \frac{\|G \1_{T}\|_{L^{1}}}{\mu(T)}
\end{align*}
so that the sizes $\|\cdot\|_{S^{1}}$, $\|\cdot\|_{S_{e}}$, and $\|\cdot\|_{S_{m}}$ satisfy the product
estimate \eqref{eq:holder-size-condition}. It follows from the
iterated Hölder inequality \eqref{eq:outer-holder}that 
\begin{align*}
\|F A\|_{L^{1}\L^{1}(S^{1})}\lesssim \|F \|_{L^{p}\L^{q}(S_{e})}\|A\|_{L^{p'}\L^{q'}(S_{m})}
\end{align*}
for conjugate exponents $\frac{1}{p}+\frac{1}{p'}=1$ and
$\frac{1}{q}+\frac{1}{q'}=1$. Furthermore we may apply
\ref{cor:iter-radon-nikodym} to $(\mb X,\nu)$ with the
local size $\|\cdot\|_{\L^{1}(S^{1})}$ so \eqref{eq:outer-duality} follows.

  \subsubsection{Interpolation}\label{sec:outer-measure-interpolation}
  Here we recall some interpolation properties of outer measure
  $L^{p}$ spaces from \cite{do2015lp} and extend them to iterated outer
  measure $L^{p}$ spaces. 

  The proof of the following Propositions can be found in
  \cite{do2015lp}.

  \begin{prop}[Logarithmic convexity of $L^{p}$ norms]\label{prop:log-convex-Lp}
    Let $(\mb X,\mu)$ be an outer measure space with size
    $\|\cdot\|_{S}$ and let $G$ be a Borel function on $\mb X$. For
    every $\theta\in(0,1)$ and for
    $\frac{1}{p_{\theta}}=\frac{1-\theta}{p_{0}}+\frac{\theta}{p_{1}}$
    with $p_{0},p_{1}\in(0,\infty]$, $p_{0}\neq p_{1}$ the inequality
  \begin{align*}
    \|G\|_{L^{p_{\theta}}(S)}\leq C_{\theta,p_{0},p_{1}}
    \|G\|_{L^{p_{0},\infty}(S)}^{1-\theta}    \|G\|_{L^{p_{1},\infty}(S)}^{\theta}
  \end{align*}
  holds.
\end{prop}

The following straight-forward remarks are useful to be able to compare outer measure
spaces with differing sizes. 

\begin{rem}[Monotonicity of outer $L^{p}$ spaces]\label{rem:Lp-monoton}
  Consider an outer measure space $(\mb X,\mu)$ with two sizes
  $\|\cdot\|_{S}$ and $\|\cdot\|_{S'}$.  Suppose that given two Borel
  functions $G$ and $G'$ on $\mb X$ we have that
  \begin{align*}
    \|G \1_{\mb X\setminus E}\|_{S}\lesssim     \|G' \1_{\mb X\setminus E}\|_{S'}
  \end{align*}
  for any $E=\bigcup_{n\in\N} T_{n}$ that is countable union of
  generating sets $T_{n}\in \mb T$. Then
  \begin{align*}
    &\|G\|_{L^{p}(S)}\lesssim \|G\|_{L^{p}(S')} 
  \end{align*}
  for all $p\in(0,\infty]$ and for iterated spaces
  \begin{align*}
    \|G\|_{L^{p}\L^{q}(S)}\lesssim \|G'\|_{L^{p}L^{q}(S')}
  \end{align*}
  for all $p,q\in(0,\infty]$. Similar statements hold for weak spaces.
\end{rem}

\begin{rem}[Interpolation of sizes]\label{rem:interp-size}
  Given an outer measure space $(\mb X,\mu)$ with two sizes
  $\|\cdot\|_{S}$ and $\|\cdot\|_{S'}$, define the sum size as 
  $\|\cdot\|_{S+S'}:=\|\cdot\|_{S}+\|\cdot\|_{S'}$. Then the following
  inequality holds for any Borel function $G$ and for any $p\in(0,\infty]$
  \begin{align*}
    &\|G\|_{L^{p}(S+S')}\leq 2\left(\|G\|_{L^{p}(S)}+\|G\|_{L^{p}(S')}\right).
  \end{align*}
\end{rem}
The proofs of the above remarks consists of simply applying the
definition of outer measure $L^{p}$ quasi-norms and as such are left to
the reader.

As a consequence of the above properties, given a function $G$ the following inequality holds:
\begin{align*}
  &\|G\|_{L^{p}\L^{q}(S)}\leq C_{q,q_{0},q_{1}}\left(\|G\|_{L^{p}\L^{q_{0},\infty}(S)}+\|G\|_{L^{p}\L^{q_{1},\infty}(S)}\right)
\end{align*}
for all $q_{0},q_{1}\in(1,\infty]$ and $q\in(q_{0},q_{1})$.

Finally we state a version of the Marcinkiewicz interpolation for maps
into outer measure $L^{p}$ spaces

  \begin{prop}[Marcinkiewicz interpolation]\label{prop:outer-marcinkiewicz}
  Let $(Y,\mc L)$ be a classical measure space, $(\mb X,\mu)$ be an
  outer measure space with size $\|\cdot\|_{S}$ and assume $1\leq p_{0}<p_{1}\leq \infty$. Let
  $T$ an operator that maps $L^{p_{0}}(Y,\mc L)+L^{p_{1}}(Y,\mc L)$ to
  Borel function on $\mb X$ so that 

\begin{description}
  \item[Scaling] $\left|T(\lambda f)\right| = \left|\lambda
      T(f)\right|$ for all $f\in L^{p_{0}}\left(Y,\mc L\right)+L^{p_{1}}\left(Y,\mc L\right)$ and
    $\lambda\in \R$;
  \item[Quasi sub-additivity]~\\
    $\left|T(f+g)\right|\leq
      C\left(\left|T(f)\right|+\left|T(g)\right|\right)$ for all
      $f,g\in L^{p_{0}}\left(Y,\mc L\right)+L^{p_{1}}\left(Y,\mc L\right)$;
    \item[Boundedness]
      \begin{equation*}
        \begin{aligned}
          &\left\|T(f)\right\|_{L^{p_{0},\infty}(S)} \leq C_{1} \left\|f\right\|_{L^{p_{0}}\left(Y,\mc L\right)}&&\forall f \in L^{p_{0}}\left(Y,\mc L\right) \\
          &\left\|T(g)\right\|_{L^{p_{1},\infty}(S)} \leq C_{2} \left\|g\right\|_{L^{p_{1}}\left(Y,\mc L\right)} &&\forall g \in L^{p_{1}}\left(Y,\mc L\right).
        \end{aligned}
      \end{equation*}
\end{description}
Then for all $f\in
L^{p_{0}}\left(Y,\mc L\right)\cap L^{p_{1}}\left(Y,\mc L\right)$ we have 
\[
  \left\| T(f) \right\|_{L^{p_{\theta}}(S)}\lesssim_{\theta,p_{0},p_{1}}
  C_{1}^{1-\theta}C_{2}^{\theta}\|f\|_{{L^{p_{\theta}}\left(Y,\mc L\right)}}
\]
with $\theta\in\left[0,\,1\right]$ and
$\frac{1}{p_{\theta}}=\frac{1-\theta}{p_{0}}+\frac{\theta}{p_{1}}$.
\end{prop}


\section{Wave-packet decomposition}\label{sec:var-carleson-embedding-reduction}

The main object of this section is to show inequality
\eqref{eq:var-duality-wavepacket} i.e. the domination of the
linearized variational Carleson operator via embedding maps.  The
following procedure follows the general scheme for obtaining
\eqref{eq:duality-wave-packet} \eqref{eq:C-LP}.

\begin{lem}\label{lem:restricted-1-represenatation}
  Consider any fixed parameters $d>b>0$, $0<d'<d-2b$, $d''>d+2b$, and
  a small enough $\epsilon>0$ appearing in properties
\eqref{eq:wave-packet-freq-support},  \eqref{eq:eq:wave-packet-geomtery-vanish}, and
  \eqref{eq:eq:wave-packet-geomtery-constant}. There exists a
  choice of truncated left and right wave packets
  $\Psi_{0,\eta,t}^{c_{-},c_{+},l}$ and
  $\Psi_{0,\eta,t}^{c_{-},c_{+},r}$ such that for all
  $c_{-}<c_{+}\in\R\cup\{+\infty\}$  the expansion
  \begin{align*}\num\label{eq:1-expansion}
    \1_{(c_{-},c_{+})}(\xi)= \iint_{\R\times\R^{+}} \big(\FT
    \Psi_{0,\eta,t}^{c_{-},c_{+},l}(\xi) + \FT \Psi_{0,\eta,t}^{c_{-},c_{+},r}(\xi) \big) \dd \eta
    \dd t
  \end{align*}
holds where the integral converges in locally uniformly for $\xi$ in $(c_{-},c_{+})$.
  
\end{lem}

\begin{proof}
  Let $\phi\in S(\R)$ and $\chi\in C^{\infty}_{c}(\R)$ be two
  non-negative functions such that for $\epsilon>0$ small enough, to
  be determined later the following holds 
  \begin{align*}\num\label{eq:phi-construction}
    \spt \FT \phi \subset B_{b}&& \spt \chi\subset B_{\epsilon}(d)\subset(b,+\infty)&&\iint_{\R\times\R^{+}} \FT \phi(\tilde t-\tilde{\eta}) \chi(\tilde{\eta})\dd
    \tilde{\eta} \frac{\dd \tilde t}{\tilde t}.
  \end{align*}
  A change of variable $\tilde \eta =t\eta$ and
  $\tilde t =\frac{t}{|\xi|}$, gives:
  \begin{align*}\num\label{eq:phi-normalization}
     \1_{(0,+\infty)}(\xi)=\iint_{\R\times\R^{+}}\FT    \phi_{\eta,t}(\xi) \chi(t \eta) \dd
    \eta\dd t 
   && \text{with } \phi_{\eta,t}(z) :=e^{i\eta z}t^{-1}\phi\big(\frac{z}{t}\big).
  \end{align*}
  Let $\gamma\in C^{\infty}_{c}\left([0,1+\epsilon)\right)$ so that
  \begin{align*}
  &\gamma(t)=1\text{ for } t\in\big[0,(1+\epsilon)^{-1}\big]
  &
  &\gamma(t)+\theta(1/t)=1 \text{ for } t\in\R^{+}.
\end{align*}
Such a function can be constructed by taking $\tilde{\gamma}$ to
satisfy the first two conditions and by setting
$  \gamma(t):=\frac{\tilde \gamma(t)}{\tilde \gamma(t)+\tilde
  \gamma(1/t)}$. Let us then set
\begin{align*}\num\label{eq:beta-def}
  &\beta(\xi) := \iint_{\R\times\R^{+}} \gamma(t') \FT
  \phi_{\eta',t'}(\xi)\chi(t'\eta')\dd \eta'\dd t'  \qquad \text{ so
  that }\\
&  \beta(t\xi)=  \iint_{\R\times\R^{+}} \gamma(t'/t) \FT
  \phi_{\eta',t'}(\xi)\chi(t'\eta')\dd \eta'\dd t' .
\end{align*}
Using \eqref{eq:phi-normalization} one obtains
  \begin{align*}
  \1_{(c_{-},c_{+})}(\xi)=\iiiint_{(\R\times\R^{+})^{2}} \FT
    \phi_{\eta,t}(\xi) \chi(t(\eta-c_{-})) \FT
    \phi_{\eta',t'}(\xi) \chi(t'(c_{+}-\eta')) \dd \eta' \dd t' \dd \eta \dd t,
  \end{align*}
  so the representation  \eqref{eq:1-expansion} holds with 
\begin{align*}
\num\label{eq:Psi-representation}
  &
    \FT \Psi_{0,\eta,t}^{c_{-},c_{+},l}(\xi):=\chi(t(\eta-c_{-}))\, \FT \phi_{\eta,t}(\xi)
    \,\beta(t(c_{+}-\xi))
  \\
  &
    \FT \Psi_{0,\eta,t}^{c_{-},c_{+},r}(\xi):=\chi(t(c_{+}-\eta))\,\FT
    \phi_{\eta,t}(\xi)\,\beta(t(\xi-c_{-})).
\end{align*}

It remains to check that $\Psi_{0,\eta,t}^{c_{-},c_{+},l}$ are left
truncated wave packets. By symmetry it will follow that
$\Psi_{0,\eta,t}^{c_{-},c_{+},r}$ is a right truncated wave
packet.  First of all \eqref{eq:wave-packet-freq-support} holds according to
 \eqref{eq:Psi-representation} since
 $\spt \FT \phi_{\eta,t}(\xi)\subset B_{bt^{-1}}(\eta)$.

 Notice that
 \begin{align*}\num\label{eq:beta-properties}
&   \spt \beta \subset\big(\frac{d-\epsilon-b}{1+\epsilon},
   +\infty\big) && \beta(\xi)=1 \text{ on } \big((d+\epsilon+b)(1+\epsilon),
   +\infty\big).
 \end{align*}
As a matter of fact 
 the integrand in \eqref{eq:beta-def} is non-zero
only if $t'(\xi-\eta')\in B_{b}$ and $t'\eta'\in B_{\epsilon}(d)$ so
$t'\xi\in B_{\epsilon+b}(d)$. This shows that
\begin{align*}
&  \xi\leq \frac{d-\epsilon-b}{1+\epsilon}\implies t' > 1+\epsilon
  \text{ or } t'<0\implies
  \gamma(t')=0 \implies \beta(\xi)=0\\
  &\xi\geq (d+\epsilon+b)(1+\epsilon)  \implies t' < (1+\epsilon)^{-1} \implies
  \gamma(t')=1 \implies \beta(\xi)= 1
\end{align*}
 where the last equality follows from \eqref{eq:phi-normalization}.

 We now check that \eqref{eq:eq:wave-packet-geomtery-vanish} holds. It
 follows from \eqref{eq:Psi-representation} that
 $\Psi_{y,\eta,t}^{c_{-},c_{+},l}(\xi)$ vanishes unless
 $\chi(t(\eta-c_{-}))\neq 0$ i.e. unless
 $t(\eta-c_{-})\in B_{\epsilon}(d)$. Also
 $\Psi_{y,\eta,t}^{c_{-},c_{+},l}(\xi)=0$ unless $t(\xi-\eta)>-b$ and
 $t(c_{+}-\xi)>\frac{d-\epsilon-b}{1+\epsilon}$ i.e. unless
 $t(c_{+}-\eta)>\frac{d-\epsilon-b}{1+\epsilon}-b$ As long as
 $0<d'<d-2b$ one can choose $\epsilon>0$ small enough for
 \eqref{eq:eq:wave-packet-geomtery-vanish} to hold.

We now check that \eqref{eq:eq:wave-packet-geomtery-constant}
holds. We have that $\beta(t(c_{+}-\xi))=1$ if
$t(c_{+}-\xi)>(d+\epsilon+b)(1+\epsilon)$ and we know that $\FT
\phi_{\eta,t}(\xi)\neq 0$ only if 
$t(\xi-\eta)\in B_{b}$ thus  if
$t(c_{+}-\eta)>(d+\epsilon+b)(1+\epsilon)+b$ then
\begin{align*}
  \FT \Psi_{0,\eta,t}^{c_{-},c_{+},l} = \chi(t(\eta-c_{-}))\, \FT \phi_{\eta,t}(\xi)
  =:   \FT \Psi_{0,\eta,t}^{c_{-},+\infty,l}
\end{align*}
so \eqref{eq:eq:wave-packet-geomtery-constant} holds as long as
$d''>d+2b$ and $\epsilon>0$ is chosen small enough

We now need to check the smoothness conditions
\eqref{eq:wave-packet-smoothness}.  We must show
that the functions
\begin{align*}
  &\FT \Psi_{0,\eta,t}^{c_{-},c_{+},l}\left(\frac{\xi+\eta}{t}\right)
  &
  &t^{-1}\partial_{c_{-}}\FT \Psi_{0,\eta,t}^{c_{-},c_{+},r}\left(\frac{\xi+\eta}{t}\right)
  &
  &t^{-1}\partial_{c_{+}}\FT \Psi_{0,\eta,t}^{c_{-},c_{+},r}\left(\frac{\xi+\eta}{t}\right)
\end{align*}
are all uniformly bounded in $S(\R)$ for all $\eta,t\in\R\times\R^{+}$
and $c_{-}<c_{+}\in \R$. Clearly 
\begin{align*}
  &\Psi_{0,\eta,t}^{c_{-},c_{+},l}\Big(\frac{\xi+\eta}{t}\Big)=
    \chi(t\eta-tc_{-}) \,\FT
    \phi_{\eta,t}\Big(\frac{\xi+\eta}{t}\Big)\,\beta(t c_{+}-\xi+\eta)
  \end{align*}
and the claim follows.
\end{proof}

\begin{cor}\label{cor:enlarged-parameter-range}
  Let us fix a set of
  parameters $d'',d',d>0$ with $d''>\max(d';d)$ and $3d>d'$. Then
  for any $\epsilon>0$ small enough there exists $b>0$ such that
  there exists a choice of left and right truncated wave packets
  $\Psi_{0,\eta,t}^{c_{-},c_{+},l}$ and $\Psi_{0,\eta,t}^{c_{-},c_{+},r}$ such
  that \eqref{eq:1-expansion} holds for all $c_{-}<c_{+}\in
  \R\cup\{+\infty\}$. 
\end{cor}
\begin{proof}
    If $d''>d>d'>0$ then let us choose $\epsilon>0$ and $b>0$ small
  enough so that the conditions for Lemma
  \ref{lem:restricted-1-represenatation} hold. Then the Lemma provides us with
  wave packets $\Psi_{0,\eta,t}^{c_{-},c_{+},l}$ and
  $\Psi_{0,\eta,t}^{c_{-},c_{+},r}$ such that \eqref{eq:1-expansion}
  holds as required.

  Suppose now that $3d>d'\geq d$ and $d''>d'$ and consider the set of
  parameters    $\tilde{d}'',\tilde{d}',\tilde{d},\tilde{b},
  \tilde{\epsilon}>0$ given by 
\begin{align*}
  &\tilde{\epsilon}=\epsilon
  &&  \tilde{b}= b-\delta 
  &&  \tilde{d} = d +\delta
  &&  \tilde{d}'=d' -\delta 
     && \tilde{d}''=d'' -\delta
\end{align*}
for some $d>\delta>0$.  We need to check that the above parameters
satisfy the assumptions of Lemma
\ref{lem:restricted-1-represenatation} that will give us the left and
right truncated wave-packets $\tilde{\Psi}_{0,\eta,t}^{c_{-},c_{+},l}$
and $\tilde{\Psi}_{0,\eta,t}^{c_{-},c_{+},r}$ for which
\eqref{eq:wave-packet-freq-support},
\eqref{eq:eq:wave-packet-geomtery-vanish}, and
\eqref{eq:eq:wave-packet-geomtery-constant} hold with these modified
parameters As long as $2\delta+\tilde{b}<\tilde{d}-\tilde{\epsilon}$,
setting
$\Psi_{0,\eta,t}^{c_{-}c_{+},l}:=\tilde{\Psi}_{0,\eta+\delta
  t^{-1},t}^{c_{-},c_{+},l}$ and
$\Psi_{0,\eta,t}^{c_{-}c_{+},r}:=\tilde{\Psi}_{0,\eta-\delta
  t^{-1},t}^{c_{-},c_{+},r}$ will provide us with the required
wave-packets so that  \eqref{eq:1-expansion} holds.

Set $b=\frac{d'-d}{2(1-3\epsilon)}$ and $\delta=(1-\epsilon)b$ so
that $\tilde{b}=\epsilon b$ with  $\epsilon>0$ small enough for the
subsequent inequalities to hold.  All the abovementioned conditions
hold since 
\begin{align*}
    &\tilde{d}-\epsilon-\tilde{b}-2\delta =d+\delta-\epsilon -b +\delta
    -2\delta=d-\epsilon-b=d-\epsilon - \frac{d'-d}{2(1-3\epsilon)}>0\\
& \tilde{b}=\epsilon b>0\\
  & \tilde{d}-\tilde{b} = d -b +2\delta>0\\
  &\tilde{d}'>d'-\delta = d'-\frac{1-\epsilon}{2}\frac{d'-d}{1-3\epsilon}>0\\
  &\tilde{d} -2\tilde{b}-\tilde{d}' = d-d'-2b +4\delta =
    d-d'+2(1-2\epsilon)b=
    (d'-d)\left(\frac{1-2\epsilon}{1-3\epsilon}-1\right)>0\\
  &\tilde{d}''-\tilde{d}-2\tilde{b} = d''-d-2b> d''-d' + (d'-d)\left(1-\frac{1}{1-3\epsilon}\right)>0.
\end{align*}
This concludes the proof. 
\end{proof}

As a consequence we obtain the following representation Lemma.

\begin{lem} \label{lem:jump-presentation} Let us fix a set of
  parameters $d'',d',d>0$ with $d''>\max(d';d)$ and $3d>d'$.  For any
  $\epsilon>0$ small enough there exists $b>0$ such that for any
  $f\in S(\R)$ and $c_{-}<c_{+}\in \R\cup \{+\infty\}$ the expansion
  \begin{align*}\num\label{eq:f-1-wave-packet-representation}
    \int_{c_{-}}^{c_{+}}\FT f(\xi) e^{i\xi x}\dd \xi = \iiint_{\mb X}
    f*\psi_{\eta,t}(y) \Big(\Psi_{y,\eta,t}^{c_{-},c_{+},l}(z) +
    \Psi_{y,\eta,t}^{c_{-},c_{+},r}(z)\Big) \dd y \dd \eta \dd t
  \end{align*}
  holds. Here $\Psi_{y,\eta,t}^{c_{-},c_{+},l}$ and
  $\Psi_{y,\eta,t}^{c_{-},c_{+},l}$ are some   
  left and right truncated wave packets for which properties
  \eqref{eq:wave-packet-smoothness},
  \eqref{eq:wave-packet-freq-support},
  \eqref{eq:eq:wave-packet-geomtery-vanish}, and
  \eqref{eq:eq:wave-packet-geomtery-constant} hold with the
  parameters above.  The function  $\psi_{\eta,t}$ is obtained from
  some $\psi\in S(\R)$ as in
  \eqref{eq:psi-eta-t}; we also have 
  \begin{align*}
    &\spt \FT \psi \in B_{(1+\epsilon)b}\qquad \text{with}  (1+\epsilon)b < d-\epsilon.
  \end{align*}
\end{lem}
\begin{proof}
 Let us choose $\psi\in S(\R)$ such that $\spt \FT \psi \in
  B_{(1+\epsilon)b}$ and $\FT \psi=1$ on $B_{b}$ so that
  \begin{align*}
&    \FT \Psi_{0,\eta,t}^{c_{-},c_{+},l}(\xi)= \FT \psi_{\eta,t}(\xi)\FT \Psi_{0,\eta,t}^{c_{-},c_{+},l}(\xi)&&    \FT \Psi_{0,\eta,t}^{c_{-},c_{+},r}(\xi)= \FT \psi_{\eta,t}(\xi)\FT \Psi_{0,\eta,t}^{c_{-},c_{+},r}(\xi)
  \end{align*}
  and let us set
  $\Psi_{y,\eta,t}^{c_{-},c_{+},l}(z)=\Psi_{0,\eta,t}^{c_{-},c_{+},l}(z-y)$
  and
  $\Psi_{y,\eta,t}^{c_{-},c_{+},r}(z)=\Psi_{0,\eta,t}^{c_{-},c_{+},r}(z-y)$. It
  follows that 
  \begin{align*}
    &    \iiint_{\mb X}
      f*\psi_{\eta,t}(y) \Big(\Psi_{y,\eta,t}^{c_{-},c_{+},l}(z) +
      \Psi_{y,\eta,t}^{c_{-},c_{+},r}(z)\Big) \dd y \dd \eta \dd t
    \\
    &=
      \iint_{\R\times\R^{+}}
      f*\psi_{\eta,t}(y) *\Big(\Psi_{0,\eta,t}^{c_{-},c_{+},l} +
      \Psi_{0,\eta,t}^{c_{-},c_{+},r}\Big)(z)  \dd \eta \dd t
    \\
    &=
      \Fourier^{-1}\Bigg(\iint_{\R\times\R} \FT f(\xi) \FT
      \psi_{\eta,t}(\xi )\Big(\FT \Psi_{0,\eta,t}^{c_{-},c_{+},l}(\xi) +
      \FT \Psi_{0,\eta,t}^{c_{-},c_{+},r}(\xi)\Big)\Bigg)=
      \Fourier^{-1}\Big(\FT f(\xi) \1_{c_{-},c_{+}}(\xi)\Big)
  \end{align*}
as required, where $\Fourier^{-1}$ is the inverse Fourier transform.
\end{proof}

As a corollary of the above Lemma we have the following pointwise wave-packet
representation for the linearized variational Carleson operator:
\begin{align*}  
  \sum_{k\in\Z}a_{k}(z)&
  \int_{\mf c_{k}(z)}^{\mf c_{k+1}(z)} \FT f(\xi) e^{i\xi z}\dd \xi
  \\=&
    \sum_{k\in\Z} \iiint_{\mb X} f*\psi_{\eta,t}(y)
    \Big(\Psi^{\mf c_{k}(z),\mf c_{k+1}(z),l}_{y,\eta,t}(z)+\Psi^{\mf c_{k}(z),\mf c_{k+1}(z),r}_{y,\eta,t}(z)\Big)a_{k}(z)\dd \eta\dd y \dd t  .
\end{align*}
Setting
\begin{align*}
  \mb A_{\mf c}(y,\eta,t):=\int_{\R}\sum_{k\in\Z} \Big(\Psi^{\mf
  c_{k}(z),\mf c_{k+1}(z),l}_{y,\eta,t}(z)+\Psi^{\mf c_{k}(z),\mf
  c_{k+1}(z),r}_{y,\eta,t}(z)\Big)a_{k}(z) \dd z
\end{align*}
gives \eqref{eq:var-duality-linearized}. We also remark that if $\c $
and $\a$ are as in \eqref{eq:var-to-carleson-reduction} then the above
construction reduces to the one described by \eqref{eq:C-LP},
\eqref{eq:F} and \eqref{eq:A} thus showing
\eqref{eq:carleson-duality-linearized}.

Finally notice that if we fix $\Theta=(\alpha^{-},\alpha^{+})=(-1,1)$
and set $d<1<d'<d''$ with $d'<3d$, then for every $\epsilon>0$ small
enough we may apply Lemma \ref{lem:jump-presentation} to obtain the
parameter $b>0$ and wave-packets
$\Psi^{\mf c_{-},c_{+},l}_{y,\eta,t}(z)$,
$\Psi^{\mf c_{-},c_{+},l}_{y,\eta,t}(z)$, and
$\psi_{\eta,t}$. Supposing that $\epsilon>0$ is small enough so that
$d+\epsilon<\alpha^{+}=1$ we can find
$(1+\epsilon)b<\beta^{+}<d-\epsilon$ and set
$\Theta^{(i)}=(\beta^{-},\beta^{+})=(-\beta^{+},\beta^{+})$. Thus
there exists a set of parameters
$\alpha^{-}<\beta^{-}<\beta^{+}<\alpha^{+}$ such that
\eqref{eq:var-duality-wavepacket} holds and
\eqref{eq:truncate-geometry} is satisfied so that Theorem
\ref{thm:max-var-mass-bounds} and Theorem \ref{thm:energy-boundedness}
hold.


\section{The auxiliary embedding map}\label{sec:aux-embedding}

In this section we introduce an auxiliary embedding map used to
control the embedded function $\mb A$.  The bounds with the same
exponents as in \eqref{eq:var-mass-embedding-bound} hold for the
auxiliary embedded function $\mb M$ with $S^{\infty}$ in lieu of
$S_{m}$. However it is technically easier to
control the super-level outer measure $\mu\left(\|\mb M\|_{S^{\infty}}>\lambda\right)$ of the
auxiliary embedded function $\mb M$. A crucial covering Lemma implies
non-iterated outer $L^{p'}$ space bounds for $\mb M$ while a locality
property and a projection Lemma allows for the extention to iterated
outer $L^{p'}\L^{q}$ spaces.

The auxiliary embedding map
associates to $\a\in C^{\infty}_{c}(l^{r'})$ the
function on $\mb X$ given by
\begin{align*}\num \label{eq:aux-map}
  \mb M  (y,\eta,t):= \int_{\R} \Big(\sum_{k\in\Z}|\a_{k}(z)|^{r'}
  \,\1_{\Theta}\left(t\left(\eta-\mf c_{k}(z)\right)\right)\Big)^{1/r'}
  W_{t}(z-y) \dd  z
\end{align*}
where the bump function $W$ is as in \eqref{eq:def:bump}. 

\begin{prop}[Bounds on the auxiliary embedding map $\mb M$]\label{prop:aux-map-bounds}
  For any $r'\in [1,\infty]$, $p'\in(1,\infty]$, and $q'\in(r',\infty]$
  and for any function $\a\in L^{p'}(l^{r'})$ the function
  $\mb M$ defined by \eqref{eq:aux-map} satisfies the bounds
\begin{align*}\num\label{eq:aux-bounds}
  &\|\mb M  \|_{L^{p'}\L^{q'}(S^{\infty})}
    \lesssim \|\a\|_{L^{p'}(l^{r'})}
\end{align*}
 where $S^{\infty}(\mb
M):=\sup_{(y,\eta,t)\in\mb X} \mb M (y,\eta,t)$. Furthermore the weak
endpoint bounds
\begin{align*}\num\label{eq:aux-bounds-endpoint}
  &\|\mb M  \|_{L^{p'}\L^{r',\infty}(S^{\infty})}\lesssim
    \|\a\|_{L^{p'}(l^{r'})} &&p'\in\left(1,\,\infty\right]\\
  &\|\mb M  \|_{L^{1,\infty}\L^{q'}(S^{\infty})}\lesssim
    \|\a\|_{L^{1}(l^{q'})}&& q'\in(r',\infty]\\
  &\|\mb M  \|_{L^{1,\infty}\L^{r',\infty}(S^{\infty})} 
    \lesssim \|\a\|_{L^{1}(l^{r'})}.
\end{align*}
hold. All the above inequalities hold as
long as $N>0$ in \eqref{eq:def:bump} is large enough and with constants independent of
the stopping sequence $\mf c$ appearing in \eqref{eq:aux-map}.
\end{prop}

We may make two reductions to prove the above bounds.  First of all
one can substitute $W_{t}(z)$ by a normalized characteristic function
of a ball. As a matter of fact set
\begin{align*}
\mb M_{R}
  (y,\eta,t):=\fint_{B_{Rt}(y)}\Big(\sum_{k\in\Z}|\a_{k}(z)|^{r'}\,\1_{\Theta}\big(t(\eta-\c_{k}(z))\big)\Big)^{1/r'} \dd  z 
\end{align*}
so that
$\mb M (y,\eta,t)\lesssim\sum_{n\in\N} R^{-Nn}\,\mb M_{R^{n}}
(y,\eta,t)$. Thus it is sufficient to prove that the bounds
\eqref{eq:aux-bounds} hold for $\mb M_{R}$ with a constant that grows
at most as $R^{N'}$ for some $N'>0$ as $R\to\infty$. The bounds for
$\mb M$ follow by quasi-subadditivity as remarked in Section
\ref{sec:outer-dominated-convergence} as long as $N>N'$. For the
second reduction split $\Theta=\Theta^{+}\cup\Theta^{-}$ into
$\Theta^{+}:=\Theta\cap [0,+\infty]$ and
$\Theta^{-}:=\Theta\cap[-\infty,0]$. Set
\begin{align*}\num\label{eq:aux-map-R-plus}
  & \mb  M^{\pm}_{R} 
    (y,\eta,t):=\fint_{B_{Rt}(y)}\Big(\sum_{k\in\Z}|\a_{k}(z)|^{r'}\,\1_{\Theta^{\pm}}\big(t(\eta-\c_{k}(z))\big)\Big)^{1/r'}
    \dd  z
\end{align*}
so that $ \mb M_{R} \leq \mb M^{+}_{R} + \mb M^{-}_{R}$.  Thus it will
suffice to provide the proof of the bounds \eqref{eq:aux-bounds} only
for $\mb M_{R}^{+}$

We begin by introducing the concept of disjoint tents relative to the
embedding \eqref{eq:aux-map-R-plus} and record an important covering
lemma.

\begin{defn}[$Q^{+}$-disjointness]
  Let $Q>0$.  We say two tents 
  $T(x,\xi,s)$ and $T(x',\xi',s')$ are $Q^{+}$-disjoint if
  either
\begin{align*}
  &B_{Qs}(x)\cap B_{Qs'}(x')=\emptyset
  &&\text{or }
  && \{c \sthat s(\xi-c)\in\Theta^{+}\}
     \cap\{c \sthat s'(\xi'-c)\in\Theta^{+}\}=\emptyset.
\end{align*}
\end{defn}
Notice that if a sequence of tents $T(x_{l},\xi_{l},s_{l})_{l\in\N}$ are pairwise
$Q^{+}$-disjoint, with $Q\geq R$, then for every $z\in\R$
\begin{align*}
  \left|\sum_{l\in\N} \1_{\Theta^{+}}\big(s_{l}(\xi_{l}-\c_{k}(z))\big)
  \1_{B_{R}}\left(\frac{x_{l}-z}{s_{l}}\right)\right|\leq 1
\end{align*}
and the bound
\begin{align*}\num\label{eq:aux-step-bound}
  &\sum_{l\in\N} s_{l}\,\mb M^{+}_{R}(x_{l},\xi_{l},s_{l})^{r'} 
  \leq\sum_{l\in\N} s_{l} \fint_{B_{Rs_{l}}(x_{l})}\sum_{k\in\Z}|\a_{k}(z)|^{r'}\,\1_{\Theta^{\pm}}\big(t(\xi_{l}-\mf
    c_{k}(z))\big)\dd z\\
  &\leq (2R)^ {-1}\int_{\R}\|\a(z)\|_{l^{r'}}^{r'} \dd z =(2R)^ {-1} \|\a\|_{L^{r'}(l^{r'})}^{r'}
      \end{align*}
      holds.

      What follows is a covering lemma. We remark that this is the only instance where we
      require smoothness and rapid decay assumptions on $\a$.
      
\begin{lem}\label{lem:superlevel-covering}
  Let $\a\in C^{\infty}_{c}(l^{r'})$. If $Q>R>R_{0}$ for some $R_{0}>0$
  depending on $\Theta$ the super level set 
  \begin{align*}
    E_{\lambda,R}:=\left\{(x,\xi,s) \sthat \mb
    M^{+}_{R}(x,\xi,s)\geq\lambda\right\}
  \end{align*}
  admits  a finite covering
  $\bigcup_{l=1}^{L} 3 Q^{2}\, T_{l} \supset E_{\lambda,R}$ with tents
  $Q^{+}$-disjoint tents  $T_{l}=T(x_{l},\xi_{l},s_{l})$ centered at points 
  $(x_{l},\xi_{l},s_{l})\in E_{\lambda,R}$.
\end{lem}

\begin{proof}
  Introduce the relation $\triangleleft$ between points of $\mb X$
  such that $(x,\xi,s)\triangleleft(x',\xi',s')$ if
  $B_{Qs}(x)\cap B_{Qs'}(x')\neq \emptyset$, $s(\xi-\xi')\in\Theta$
  and $s'>Qs$.  We say $(x,\xi,s)$ is maximal in a set $P\subset\mb X$
  if there is no $(x',\xi',s')\in P$ such that
  $(x,\xi,s)\triangleleft (x',\xi',s')$. Notice that $E_{\lambda,R}$
  is $(x,t)$-bounded in the sense that for some $C>1$ large enough
\begin{align*}
  E_{\lambda,R}\subset B_{C}(0)\times \R\times (0,C)
\end{align*}
holds. As a matter of fact $\mb M_{R}^{+}(y,\eta,t)\lesssim
(Rt)^{-1} \|\a\|_{L^{1}(l^{r'})}$ and 
$\mb M_{R}^{+}(y,\eta,t)=0$ if
$\dist(y;\,\spt \a)>tR$ so if $(y,\eta,t)\in E_{\lambda,R}$ then
$t<C$ and $|y|<C$ for some $C>0$ depending on
$\a$. Thus any non-empty subset $P\subset E_{\lambda,R}$
admits a maximal element. 
      
Inductively construct a covering starting with an empty collection of
tents $\mc T^{0}=\emptyset$. At the $l$\textsuperscript{th} step
consider the points in the set 
\begin{align*}\num\label{eq:aux-selection-condition}
 E_{\lambda,R}\setminus\bigcup_{T\in\mc T^{l-1}} 3 Q^{2}\, T
\end{align*}
and select from it a point $(x_{l},\xi_{l},s_{l})$ that is maximal with respect
to the relation $\triangleleft$ and set 
and $\mc T^{l} = \mc T^{l-1}\cup \{T(x_{l},\xi_{l},s_{l})\}$. We claim that at each
step of the algorithm all the selected tent $T(x_{l},\xi_{l},s_{l})$
are pairwise $Q^{+}$-disjoint.  Reasoning by contradiction, suppose that
two tents $T(x_{l},\xi_{l},s_{l})$ and $T(x_{l'},\xi_{l'},s_{l'})$ with
$l<l'$ are not $Q^{+}$-disjoint, then
$B_{Qs_{l}}(x_{l})\cap B_{Qs_{l'}}(x_{l'})\neq \emptyset$ and there
also exists a $c\in\R$ such that $s_{l}(\xi_{l} - c) \in \Theta^{+}$
and $s_{l'}(\xi_{l'}-c)\in\Theta^{+}$. Recall that $\Theta^{+}=[0,\alpha^{+}]$ so 
\begin{align*}
&\xi_{l}-s_{l}^{-1}\alpha^{+}\leq c\leq \xi_{l}
  &&\xi_{l'}-s_{l'}^{-1}\alpha^{+}\leq c\leq \xi_{l'}
\end{align*}
If $s_{l'}\geq Q s_{l}$  one would have
\begin{align*}
  -s_{l}^{-1}Q^{-1}\alpha^{+}\leq -s_{l'}^{-1}\alpha^{+}\leq\xi_{l}-\xi_{l'}\leq s_{l}^{-1}\alpha^{+}
\end{align*}
and thus $s_{l}(\xi_{l}-\xi_{l'})\in \Theta$ as long as $\alpha^{-}\leq-R_{0}^{-1}\alpha^{+}$. This contradicts the maximality
of $(x_{l},\xi_{l},s_{l})$  that was chosen before
$(x_{l'},\xi_{l'},s_{l'})$. On the other hand if $s_{l'}<Q s_{l}$ then
\begin{align*}
  - s_{l}^{-1}\alpha^{+}\leq \xi_{l'}-\xi_{l}\leq s_{l'}^{-1}\alpha^{+}
\end{align*}
and, as long as $Q\geq R_{0}\geq \alpha^{+}$, this implies that
$(x_{l'},\xi_{l'},s_{l'})\in 3Q^{2} T(x_{l},\xi_{l},s_{l})$
contradicting the selection condition.

Finally notice that the selection algorithm terminates after
finitely many steps since at every step
\eqref{eq:aux-step-bound} holds having chosen $Q\geq R$, since $s_{l}$
are bounded from below since 
$\mb M^{+}_{R}(x_{l},\xi_{l},s_{l})\geq \lambda$.  Thus
$E_{\lambda}\subset \bigcup_{l=1}^{L}3 Q^{2}T_{l}$.
\end{proof}

A consequence of the above Lemma are non-iterated bounds for $\mb M^{+}_{R}$.
\begin{prop}\label{prop:aux-map-non-iter}
  Given $\a\in L^{p'}(l^{r'})$ with $p'\in(r',\infty]$ the bound
  \begin{align*}\num\label{eq:aux-bounds-non-iter}
          &\|\mb M^{+}_{R}\|_{L^{p'}(S^{\infty})}
            \lesssim_{R} \|\a\|_{L^{p'}(l^{r'})}
  \end{align*}
  holds. Furthermore the weak endpoint bound
  \begin{align*}\num\label{eq:aux-bounds-non-iter-weak}
          &\|\mb M^{+}_{R}\|_{L^{r',\infty}(S^{\infty})}\lesssim_{R}\|\a\|_{L^{r'}(l^{r'})}
  \end{align*}
  holds. All the above bounds hold with a constant that grows at most
  polynomially in $R$ as $R\to \infty$ and is independent of the
  stopping sequence $\mf c$ appearing in \eqref{eq:aux-map-R-plus}.
\end{prop}
      
The bound  \eqref{eq:aux-bounds-non-iter} for $p=\infty$ is straightforward:
\begin{align*}
  & \mb M^{+}_{R}  (y,\eta,t)= \fint_{B_{tR}(y)} \Big(\sum_{k\in\Z}|\a_{k}(z)|^{r'}
  \,\1_{\Theta^{+}}\big(t(\eta-c_{k}(z))\big)\Big)^{1/r'}\dd  z 
 \\
 &
 \leq \int_{\R}\sum_{k\in\Z}|\a_{k}(z)|^{r'}
   \,t^{-1}\1_{B_{R}}\left(\frac{z-y}{t}\right)\dd 
  z  \lesssim  \|\a\|_{L^{\infty}(l^{r'})}^{r'}.
\end{align*}
It is sufficient to show bound \eqref{eq:aux-bounds-non-iter-weak} so
that will \eqref{eq:aux-bounds-non-iter} follow for $p\in(r',\infty)$ by
interpolation \ref{prop:outer-marcinkiewicz}.  In particular to obtain
\eqref{eq:aux-bounds-non-iter-weak} we will show that given
$\lambda>0$ the bound on the measure of the super-level set
\begin{align*}\num\label{eq:M-om-R-r'-bound}
& \mu(E_{\lambda,R})\lesssim_{R} \lambda^{-r'}\|\a\|_{L^{r'}(l^{r'})}^{r'}
\end{align*}
holds. It is sufficient to consider the covering provided by Lemma \ref{lem:superlevel-covering}
with $Q=R$.  Since $(x_{l},\xi_{l},s_{l})\in E_{\lambda,R}$ and the
covering $\mc T=\mc T^{L}=\{T(x_{l},\xi_{l},s_{l})\}_{l\in L}$ consists of $Q^{+}$-disjoint tents, the bound
\eqref{eq:aux-step-bound} gives
\begin{align*}
  \lambda^{r'} \sum_{l=1}^{L}s_{l} \leq (2R)^{-1}\|\a\|_{L^{r'}(l^{r'})}^{r'}.
\end{align*}
Since $E_{\lambda,R}\subset
\bigcup_{l=1}^{L}3R^{3}T(x_{l},\xi_{l},s_{l})$ one deduces
\begin{align*}
  \mu(E_{\lambda,R}) \lesssim_{R}
  \sum_{l=1}^{L}\mu\big(T(\xi_{l},x_{l},s_{l})\big) \leq
  \sum_{l=1}^{L}s_{l}
  \leq  \frac{\|\a\|_{L^{r'}(l^{r'})}^{r'}}{\lambda^{r'}}
\end{align*}
where the implied constant grows polynomially in $R$ as required.

The proof of \ref{prop:aux-map-bounds} relies on a locality property
and a strip projection lemma.

\begin{lem}[Locality of $\mb M^{+}_{R}$]\label{lem:locality-M-Om-R}
  Consider a strip $D=D(x,s)$ and a function $\a\in L^{1}_{loc}(l^{r'})$ with
  \begin{align*}
    \dist\big(\spt \a;B_{s}(x)\big) > Rs
  \end{align*}
  then for all $(y,\eta,t)\in D(x,s)$ we have
  \begin{align*}
        \mb M^{+}_{R}\,\1_{D(x,s)} = 0.
  \end{align*}
\end{lem}
\begin{proof}
The statement follows directly from the definition
\eqref{eq:aux-map-R-plus} of the embedding. As a matter of fact if
$(y,\eta,t)\in D_{x,s}$ then $B_{tR}(y)\subset B_{sR}(y)$ and $B_{sR}(y)\cap \spt
\a=\emptyset $.
\end{proof}

\begin{lem}[Mass projection for $\mb M^{+}_{R}$]\label{lem:mass-projection}
  Fix any  collection of pairwise disjoint strips
  $D(\zeta_{m},\tau_{m})$, $m\in\{1,\dots,M\}$ and any finite
  collection of $Q^{+}$-disjoint tents
  \begin{align*}
    T(x_{l},\xi_{l},s_{l})\not\subset \bigcup_{m=1}^{M}
  D(\zeta_{m},3\tau_{m}),\qquad l\in\{1,\dots,L\}
  \end{align*}
  with $Q>2R>2$. Given a function $\a\in L^{1}_{loc}(l^{r'})$ and a
  stopping sequence $\c$ there exists a function $\tilde{\a}\in L^{1}_{loc}(l^{r'})$ and a
  new stopping sequence $\tilde{\c}$ such that
  \begin{align*}\num\label{eq:tilde-a-bounds}
    &\|\tilde{\a}(z)\|_{l^{r'}} \lesssim
      \fint_{B_{\tau_{m}}(\zeta_{m})}\|
      \a(z)\|_{l^{r'}}\dd z
    & &
        \forall z\in
        B_{\tau_{m}}(\zeta_{m})\qquad \forall
        m\in\{1,\dots,M\}
    \\
    &\tilde{\a}_{k}(z) =  \a_{k}(z)
    && \forall z \notin \bigcup_{m=1}^{M}   B_{\tau_{m}}(\zeta_{m})\\
       \shortintertext{and} \num\label{eq:tilde-M-bound}
    & \tilde {\mb M}{}^{+}_{2R}(x_{l},\xi_{l},s_{l})      \geq \mb
    M^{+}_{R}(x_{l},\xi_{l},s_{l})&&\forall l\in\{1,\dots,L\}.
      \end{align*}
      where $\tilde{\mb M}{}^{+}_{2R}$ is the embedded function as given
      by expression \eqref{eq:aux-map-R-plus} associated to
      $ \tilde{\a}$ with the stopping sequence $\tilde{\c}$ .
    \end{lem}

\begin{proof}
  Let us order the tents $T(x_{l},\xi_{l},s_{l})$ so that
  $\xi_{l}\leq \xi_{l'}$ if $l<l'$ . For every strip
  $D(\zeta_{m},\tau_{m})$ let
  \begin{align*}
    &\mb L_{m} :=\Big\{l\in\{1,\dots,L\} \sthat D(x_{l},Rs_{l})\cap
      D(\zeta_{m},\tau_{m})\neq \emptyset\Big\} .
  \end{align*}
Set 
\begin{align*}    &\tilde{\a}_{k}(z) = \left\{\begin{aligned}
            &\a_{k}(z) &&\text{if } z\notin \bigcup_{m}B_{\tau_{m}}(\zeta_{m})\\
            &\fint_{B_{\tau_{m}}(\zeta_{m})}\!\!\Big(\sum_{j\in\Z}|\a_{j}(z)|^{r'}
              \1_{\Theta^{+}}\big(s_{l}(\xi_{k}-\mf c_{j}(z))\big)\Big)^{1/r'}
              \!\!\!\dd z &&
              \begin{aligned}
                &\text{if }z\in B_{\tau_{m}}(\zeta_{m}) \\&\text{and }
                k\in\mb L_{m}
              \end{aligned}
              \\
              &0 &&
              \begin{aligned}
                &\text{if }z\in B_{\tau_{m}}(\zeta_{m})\\&\text{and }
                k\notin\mb L_{m}
              \end{aligned}
          \end{aligned}
        \right.\\
        &\tilde{\mf c}_{k}(z) = \left\{
          \begin{aligned}
            &\mf c_{k}(z)&& \text{if } z\notin \bigcup_{m}B_{\tau_{m}}(\zeta_{m})\\ 
            &\xi_{k}&& \text{if }z\in B_{\tau_{m}}(\zeta_{m}) \text{ and}&&
            k\in \{1,\dots,L\}\\
            &\xi_{1} && &&k<1 \\ &\xi_{L} &&&&k>L.
          \end{aligned}
                         \right.
\end{align*}
The expressions above are well defined since $D(\zeta_{m},\tau_{m})$ are
pairwise disjoint. 

The bound \eqref{eq:tilde-a-bounds} follows by the Minkowski
inequality. For $z\in B_{\tau_{m}}(\zeta_{m})$ one has 
\begin{align*}
  &\|\tilde{\a}(z)\|_{l^{r'}} = \Bigg(\sum_{k\in\Z}\Bigg(
    \fint_{B_{\tau_{m}}(\zeta_{m})} \Big(\sum_{j\in\Z} |\a_{j}(z)|^ {r'}\1_{\Theta^{+}}\big(s_{l}(\xi_{k}-\mf c_{j}(z))\big)
     \Big)^{1/r'} \dd
    z\Bigg)^{r'}\Bigg)^ {1/r'}
  \\
  &\lesssim      \fint_{B_{\tau_{m}}(\zeta_{m})} \Big(
\sum_{k\in \mb L_{m}}     \sum_{j\in\Z}
|\a_{j}(z)|^ {r'} \1_{\Theta^{+}}\big(s_{l}(\xi_{k}-\mf c_{j}(z))\big)
    \Big)^{1/r'}\dd z\leq \fint_{B_{\tau_{m}}(\zeta_{m})}\|\a(z)\|_{l^{r'}},
\end{align*}
where the last inequality holds since the tents $T(x_{l},\xi_{l},s_{l})$
are $Q^{+}$-disjoint.

It remains to show \eqref{eq:tilde-M-bound}.
Since $T(x_{l},\xi_{l},s_{l})\not\subset D(\zeta_{m},3\tau_{m})$ for any
$m$ we have that
\begin{align*}
  &B_{Rs_{l}}(x_{l})\cap B_{\tau_{m}}(\zeta_{m})\neq\emptyset \implies
    D(\zeta_{m},\tau_{m})\subset D(x_{l},2Rs_{l})\\ \shortintertext{ so set} &\mf
M_{l}=\Big\{m \sthat D(\zeta_{m},\tau_{m})\subset D(x_{l},2Rs_{l})  \Big\}.
\end{align*}
Using the definitions of $\tilde{\a}$ and $\tilde{\c}$ we obtain
      \begin{align*}
        &\tilde {\mb M}{}^{+}_{2R}(x_{l},\xi_{l},s_{l})=
          \fint_{B_{2Rs_{l}}(x_{l})} \Big(\sum_{k\in\Z}
          |\tilde{\a}_{k}(z)|^{r'}\1_{\Theta^{+}}\big(s_{l}(\xi_{l}-\tilde{\mf
          c}_ {k}(z))\big) \Big)^{1/r'}\dd z
        \\&\gtrsim(4Rs_{l})^{-1}
             \int_{B_{Rs_{l}}(x_{l})\setminus\bigcup_{m}\!\! B_{\tau_{m}}(\zeta_{m})}
            \Big(\sum_{k\in\Z}|\a_{k}(z)|^{r'}
            \1_{\Theta^{+}}\big(s_{l}(\xi_{l}-\mf c_{k}(z))\big)\Big)^{1/r'} \dd z
        \\&
            +              (4Rs_{l})^{-1}\sum_{m\in \mf M_{l}} \int_{B_{\tau_{m}}(\zeta_{m})}
            \Big(
            \sum_{k\in\mb L_{m} }
            |
            \tilde{\a}_{k}(z)|^{r'}
            \1_{\Theta^{+}}\big(s_{l}(\xi_{l}-\xi_{k})\big)\Big)^{1/r'}
            \dd z.
      \end{align*}
    Using the fact that   $T(x_{l},\xi_{l},s_{l})$ are $Q^{+}$-disjoint
    with $Q>2R$  we obtain that $s_{l}(\xi_{l}-\xi_{k})\in\Theta^{+}$,
    $z\in B_{\tau_{m}}(\zeta_{m})$, and 
    $\tilde{\a}_{k}(z)\neq 0$ only if
    $k=l$; thus 
      \begin{align*}
        & \sum_{m\in \mf M_{l}}\int_{B_{\tau_{m}}(\zeta_{m})} \Big(
          \sum_{k\in\mb L_{m} }
          |
          \tilde{\a}_{k}(z)|^{r'}
          \1_{\Theta^{+}}\big(s_{l}(\xi_{l}-\xi_{k})\big)\Big)^{1/r'}
          \dd z
          =
          \sum_{m\in \mf M_{l}}\int_{B_{\tau_{m}}(\zeta_{m})}  
          \tilde{\a}_{l}(z)         \dd z
        \\
        &
          =
          \sum_{m\in \mf M_{l}}\int_{B_{\tau_{m}}(\zeta_{m})}  \Big(\sum_{j\in\Z}|\a_{j}(z)|^{r'}
          \1_{\Theta^{+}}\big(s_{l}(\xi_{l}-\mf
          c_{j}(z))\big)\Big)^{1/r'}
          \dd z.
      \end{align*}
      This allows us to conclude that
      \begin{align*}
        \tilde {\mb M}{}^{+}_{2R}(x_{l},\xi_{l},s_{l})\gtrsim \fint_{B_{R s_{l}}(x_{l})} \Big(\sum_{k\in\Z}|\a_{k}(z)|^{r'}
            \1_{\Theta^{+}}\big(s_{l}(\xi_{l}-\mf c_{k}(z))\big)\Big)^{1/r'}
        \!\!\!\!\dd z = \mb M^{+}_{R}(x_{l},\xi_{l},s_{l}).
      \end{align*}
    \end{proof}

    We now have all the tools to prove \eqref{eq:aux-bounds} for
    $\mb M_{R}^{+}$. We proceed by interpolation, as described in
    \ref{sec:outer-measure-interpolation}, between the four (weak)
    endpoints
    \begin{align*}
      (p',q')\in\left\{(\infty,\infty),\,(\infty,r'),\,(1,\infty),\,(1,r')\right\}.
    \end{align*}
    \begin{proof}[Proof of bounds \eqref{eq:aux-bounds} for
      $\mb M^{+}_{R}$]~

      The bound for $(p',q')=(\infty,\infty)$ follows directly from
      \eqref{eq:aux-bounds-non-iter} with $p'=\infty$.

      The bound for $(p',q')=(\infty,r')$ follows from the locality
      property \ref{lem:locality-M-Om-R}. We must show that for any
      strip $D(x,s)$ one has
      \begin{align*}
        \|\mb M^{+}_{R} \1_{D(x,s)}\|_{L^{r',\infty}(S^{\infty})}^{r'}\lesssim_{R} \nu(D(x,s)) \|\a\|_{L^{\infty}(l^{r'})}^{r'}
      \end{align*}
      but due to locality and
      \eqref{eq:aux-bounds-non-iter-weak} we have that
      \begin{align*}
        \|\mb M^{+}_{R} \1_{D(x,s)}\|_{L^{r',\infty}(S^{\infty})}^{r'}\lesssim_{R}
        \|\a\, \1_{B_{2sR}(x)}\|_{L^{r'}(l^{r'})}^{r'} \leq_{R} s \|\a\|_{L^{\infty}(l^{r'})}^{r'}
      \end{align*}
      as required.
      
      The bound for $(p',q')=(1,\infty)$ makes use of the Mass
      Projection Lemma \ref{lem:mass-projection}. We need to show that
      for every $\omega>0$ there exists $K_{\omega}\subset\mb X$ 
      such that 
    \begin{align*}
      &\nu\left(K_{\omega}\right)\lesssim_{R}\omega^{-1}\|\a\|_{L^{1}(l^{r'})}&&\left\|\mb M^{+}_{R} \,\1_{\mb X\setminus K_{\omega}}\,\1_{D(x,s)}
      \right\|_{L^{\infty}(S^{\infty})} \lesssim \omega.
    \end{align*}
    for any strip $D(x,s)$.
    Let
$      \mc K_{\omega} = \{z\in\R \sthat
      M(\|\a\|_{l^{r'}})(z)>\omega\}
$ 
    where $M$ is the Hardy-Littlewood Maximal function \eqref{eq:HL-maximal-function}.
    The set $\mc K_{\tau}$ is open and in particular is a finite union
    of intervals
    $\mc K_{\omega}=\bigcup_{m=1}^{M}B_{\tau_{m}}(\zeta_{m})$. Let
    \begin{align*}
      &K_{\omega}:=\bigcup_{m=1}^{M}D(\zeta_{m},9\tau_{m}) \implies
        \nu(K_{\omega}) \lesssim \sum_{m=1}^{M} 2\tau_{m} =|\mc K_{\omega}| \lesssim \omega^{-1}\|\a\|_{L^{1}(l^{r'})}
    \end{align*}
    by the weak $L^{1}$ bound on the Hardy-Littlewood maximal function.

    For any tent $T(y,\eta,t)\not\subset D(\zeta_{m},3\tau_{m})$ apply
    Lemma \ref{lem:mass-projection} with respect to the the strips
    $\big(D(\zeta_{m},3 \tau_{m})\big)_{m\in\{1,\dots,M\}}$ and the
    one tent $T(\xi,x,s)$. By construction we obtain a function $\tilde{\a}$ such
    that $\|\tilde{\a}\|_{L^{\infty}(l^{r'})}\lesssim \omega$. Using the
    statement of the Lemma and  bound
    \eqref{eq:aux-bounds-non-iter} we have
    \begin{align*}
      \mb M^{+}_{R} (y,\eta,t) \leq \tilde{\mb M}{}^{+}_{2R}(y,\eta,t)
      \lesssim_{R} \| \tilde{\a}\|_{L^{\infty}(l^{r'})}\lesssim \omega
    \end{align*}
    as required.
    
    The proof of the case $(p',q')=(1,r')$ goes along the same
    lines. Let us suppose, without loss of generality, that
    $\a\in C^{\infty}_{c}(l^{r'}) $We need to show that for every $\omega>0$ there exists
    $K_{\omega}\subset \mb X$ such that 
    \begin{align*}
      &\nu\left(K_{\omega}\right)\lesssim_{R}\omega^{-1}\|\a\|_{L^{1}(l^{r'})}
      &
      &\left\|\mb M^{+}_{R} \,\1_{\mb X\setminus K_{\omega}}\,\1_{D(x,s)}
        \right\|_{L^{r',\infty}(S^{\infty})}^{r'} \lesssim
        \nu\big(D(x,s)\big) \omega^{r'}.
    \end{align*}
    Choose $K_{\omega}=\bigcup_{m=1}^{M}D(\zeta_{m},9 \tau_{m})$ as before.
    Let $\lambda>0$ and set
    \begin{align*}
      &\overline{E}_{\lambda,R}= E_{\lambda,R}\cap (D(x,s)\setminus
        K_{\omega}) && E_{\lambda,R}= \left\{(y,\eta,t) \sthat \mb M_{R}^{+}>\lambda\right\}.
    \end{align*}
    The Covering Lemma \ref{lem:superlevel-covering} can be applied to
    $\overline{E_{\lambda,R}}$ with $Q>2R$ sufficiently large yielding
    a covering $\big(T(x_{l},\xi_{l},s_{l})\big)_{l\in\{1,\dots,L\}}$
    such that
    $ \bigcup_{l=1}^{L}3Q^{3}T(x_{l},\xi_{l},s_{l})\supset
    \overline{E_{\lambda,R}} $ with the tents $T(x_{l},\xi_{l},s_{l})$
    that are pairwise $Q^{+}$-disjoint. Now apply the Mass Projection Lemma
    \ref{lem:mass-projection} with respect to the strips
    $\big(D(\zeta_{m},3\sigma_{m})\big)_{m\in\{1,\dots,M\}}$ and the
    tents $T(\xi_{l},x_{l},s_{l})_{l\in\{1,\dots,L\}}$. The resulting
    $\tilde{\a}$ satisfies 
    $\|\tilde{\a}\|_{L^{\infty}(l^{r'})}\lesssim \omega$ while
    \begin{align*}
      \tilde{\mb M}{}^{+}_{2R}(x_{l},\xi_{l},s_{l})\geq   \mb
      M^{+}_{R}(x_{l},\xi_{l},s_{l})\geq \lambda.
    \end{align*}
    Using the bound
    \eqref{eq:aux-step-bound} and the locality property
    \ref{lem:locality-M-Om-R} we have that
    \begin{align*}
      &      \mu\left(\overline{E_{\lambda,R}}\right) \lesssim_{R}
        \sum_{l=1}^{L}s_{l} \lesssim \lambda^{-r'}\|\tilde{\a}
        \1_{B_{2sR}(x)}(z)\|_{L^{r'}(l^{r'})}^{r'}\lesssim s \,\omega^{r'}\lambda^{-r'}.
    \end{align*}
    This concludes the proof.
    \end{proof}


\section{Proof of Theorem \ref{thm:max-var-mass-bounds}} \label{sec:main-prop-proof}

In the previous section the bounds \eqref{eq:aux-bounds} were shown to
hold for the auxiliary embedding $\mb M$. To prove  Theorem
\ref{thm:max-var-mass-bounds} it is sufficient to show that the values
of $\mb M$ control $\|\cdot\|_{S_{m}}$.  More specifically we require
the following result.

\begin{prop}\label{prop:size-control}
  Given any union of strips $K$ and a union of tents $E$ such
  that
  \begin{align*}\num\label{eq:M-lambda-bound}
    &\mb M(y,\eta,t) \leq \lambda &&\forall (y,\eta,t)\in \mb
                                             X\setminus (K\cup E)
  \end{align*}
  the bound $ \|\mb A\,\1_{\mb X\setminus (K \cup
      E)}\|_{ S_{m}}\lesssim \lambda$ holds.
\end{prop}

Assuming that the above statement holds, Theorem
\ref{thm:max-var-mass-bounds} follows by the monotonicity property of
outer $L^{p}$ sizes \ref{rem:Lp-monoton}.

The above proposition follows from showing that the
required bound holds for all local sizes:
$\|\mb A \1_{\mb X\setminus (K\cup E)}\|_{S_{m}(T)}\lesssim \lambda$. The proof 
is divided into two parts relative to showing  $L^{1}$-type bounds over
$T^{(i)}$ and $L^{(2)}$ type bounds over $T^{(e)}$ (see \eqref{eq:def:size-m}).

The former part uses crucial disjointness properties related to the
conditions \eqref{eq:eq:wave-packet-geomtery-vanish} on the truncated
wave packets.

The latter part depends on the fact that the sizes over a single tent
$T$ resembles an $L^{2}$ estimate for variational truncation of the
Hilbert transform or of a square function in the spirit of
\cite{jones2008strong}. We will elaborate on this variational estimate
in Lemma \ref{lem:var-trunc} in the following part on technical
preliminaries.

The proof also involves a crucial stopping time argument. Similarly to
the rest of the paper we avoid discretization and formulate a
continuous version of this argument that we isolate Lemma
\ref{lem:continuous-stopping-time} below.

\subsection{Technical preliminaries}

The following variational truncation bounds are a slightly modified
version of the results appearing in
\cite{jones2008strong}.

\begin{lem}[Variational truncations of singular integral operators \cite{jones2008strong}]\label{lem:var-trunc}
  For any function $H\in L^{p}(\R)$ and $\sigma\in[0,\infty)$ let us define
  the variational truncation operator
  \begin{align*}\num\label{eq:var-trunc-op}
    &\mc V^{r}_{\sigma} H (z) = \sup_{\sigma<t_{1}<\dots<t_{k}<\dots}
    \Bigg(\sum_{k}
      |H*\Upsilon_{t_{k+1}}(z)-H*\Upsilon_{t_{k}}(z)|^{r}\Bigg)^{1/r}\\
    \shortintertext{where}
    &\Upsilon\in S(\R),\; \int_{\R} \Upsilon(z)\dd z=1 \text{ and }\Upsilon_{t}(z):=t^{-1}\Upsilon\left(\frac{z}{t}\right).
  \end{align*}
  If $r>2$ and for any $p\in(1,\infty)$, above operator satisfies the bounds
  \begin{align*}\num\label{eq:var-trunc-bounds}
    &\|\mc V^{r}_{\sigma}H\|_{L^{p}}\lesssim_{r,p}    \|H\|_{L^{p}}
  \end{align*}
  and if $\sigma>0$ then
  \begin{align*}\num\label{eq:var-trunc-qausi-const}
    &\mc V^{r}_{\sigma}H(z) \lesssim_{r,p} \fint_{B_{\sigma}(z)}M\left(\mc
      V^{r}_{\sigma}H \right)(z')\dd z'\\
  \end{align*}
  where $M$ is the Hardy-Littlewood maximal function. The implicit
  constants are allowed to depend on $\Upsilon$.
\end{lem}

We record some useful properties of so called \emph{convex regions} of
tents.

\begin{defn}[Convex regions]\label{def:convex-tent}
  A convex region of a tent is a subset $\Omega \subset T(x,\xi,s)$ of a
  tent of the form 
  \begin{align*}\num\label{eq:convex-region}
    \Omega:=\bigcup_{\theta\in\Theta}\Omega_{\theta}:=\left\{(y,\xi+\theta t^{-1},t) \in T(x,\xi,s) \sthat t>\sigma_{\theta}(y)\right\}.
  \end{align*}
for some function $\sigma_{\theta}(y):\Theta\times B_{s}(x)\to [0,s]$.
\end{defn}

Given any tent $T\in \mb T$, any collection of strips $\mc D$, and any
collection of tents $\mc T$, the set
\begin{align*}
  \Omega = T\setminus \left( \bigcup_{D\in \mc D} D  \;\cup\;
  \bigcup_{T\in \mc T} T\right)\subset T
\end{align*}
is a convex region of the tent $T$. With the next lemma we show that
the bound \eqref{eq:M-lambda-bound} on a convex regions can be
extended to larger regions with scale bound $\sigma$ that is Lipschitz  in the
space variable.

\begin{lem}[Lipschitz convex regions]\label{lem:convex-region-extention}
  Let $T(x,\xi,s)\in\mb T$ be a tent and
  $\Omega=\bigcup_{\theta\in\Theta}\Omega_{\theta}\subset T(x,\xi,s)$
  be a convex region as defined in \eqref{eq:convex-region} and let us
  fix a constant $L>2$. For every $\theta\in\Theta$ such that the
  bound
  \begin{align*}
    &\mb M(y,\eta,t) \leq \lambda \qquad \forall (y,\eta,t)\in \Omega_{\theta}\\                                   
    & \text{holds for } \Omega_{\theta}=\left\{(y,\xi+\theta t^{-1},t) \in T(x,\xi,s) \sthat t>\sigma_{\theta}(y)\right\}\neq\emptyset
  \end{align*}
  there exists a Lipschitz function
  $\tilde{\sigma}_{\theta}:\R \to \R^{+}$ with Lipschitz constant
  $L^{-1}<1/2$ such that
  \begin{align*}\num\label{eq:tilde-sigma-condition-1}
    \min\big(2s;\,L^{-1}\dist(y;B_{s}(x)) \big)\leq
    \tilde{\sigma}_{\theta}(y)&\leq 2s
    && \forall y\in\R\\\num\label{eq:tilde-sigma-condition-2}
    \tilde{\sigma}_{\theta}(y)&\leq   \sigma_{\theta}(y)
    &&\forall y\in B_{s}(x)
  \end{align*}
  and
   \begin{align*}\num\label{eq:tilde-sigma-aux-bounds}
     &sW_{s}(x-y) \mb M(y,\xi+\theta t^{-1},t) \lesssim_{L} \lambda 
     &&\forall
        y\in\R,\, t\in(\tilde{\sigma}_ {\theta}(y),3s).
   \end{align*}
\end{lem}
\begin{proof}
  Fix $\theta\in\Theta$ such that $\Omega_{\theta}$ is non-empty and
  let us drop the dependence on $\theta$ from the notation by simply
  writing $\sigma(y)$ in place of $\sigma_ {\theta}(y)$. Let us set
  \begin{align*}\num\label{eq:tilde-sigma-def}
    &\tilde{\sigma}(y):=\min\big(2s;\,\tilde{\tilde{\sigma}}(y)\big)
    &&\text{with } \tilde{\tilde{\sigma}}(y):=      \inf_{y'\in B_{s}(x)}
      \max\Big(\sigma(y');\,\frac{|y-y'|}{L}\Big)
  \end{align*}
  Clearly, this defines a function on $\R$ such that conditions
  \eqref{eq:tilde-sigma-condition-1} and
  \eqref{eq:tilde-sigma-condition-2} hold. The defined function is
  $L^{-1}$-Lipschitz. It is sufficient to show that
  $\tilde{\tilde{\sigma}}$ is $L^{-1}$-Lipschitz: for any $y\in\R$ and
  $\epsilon>0$ there exits $y'\in B_{s}(x)$ such that
  \begin{align*}
    \tilde{\tilde{\sigma}}(y)\geq (1+\epsilon)^{-1} \max\left(\sigma(y');\,\frac{|y-y'|}{L}\right)
  \end{align*}
  and thus for any  $y''\in \R$ one has 
  \begin{align*}
    \tilde{\tilde{\sigma}}(y'')&\leq \max
    \Big(\sigma(y');\,\frac{|y''-y'|}{L}\Big)\leq \max
    \Big(\sigma(y');\,\frac{|y-y'|}{L}\Big) +
  \frac{|y''-y|}{L}\\
    &\leq (1+\epsilon) \tilde{\tilde{\sigma}}(y) +   \frac{|y''-y|}{L}.
  \end{align*}
  Since $\epsilon>0$ was arbitrary and one can invert the role of $y''$
  and $y$ in the above reasoning we obtain that
  $|\tilde{\tilde{\sigma}}(y'')-\tilde{\tilde{\sigma}}(y)|\leq
  \frac{|y''-y|}{L}$ as required. 

  Let us now check that \eqref{eq:tilde-sigma-aux-bounds} holds.
  Suppose that $y\in\R$ and
  $t\in(\tilde{\sigma}(y),3s]$.  Let us distinguish the cases
  $t\in(\tilde{\sigma}(y),2s)$
  and $t\in[2s,3s)$. In the first case there exists $y'\in B_{s}(x)$ and
  $t'\in(\sigma_{\theta}(y'),s)$ such that $t'\in(t/2,t)$ and
  $|y-y'|<2Lt$ and thus it follows that $|x-y|<2Ls$. It follows that 
    \begin{align*}
      &    W_{t}(z-y)\lesssim_{L} W_{t}(z-y')\lesssim W_{t'}(z-y')
      &&\forall z\in \R
      \\
      \text{thus }\quad& sW_{s}(x-y)W_{t}(z-y)\lesssim_{L}W_{t'}(z-y')&&\forall z\in \R
  \end{align*}
  In the case that $t\in[2s,3s]$ there also exists $y'\in B_{s}(x)$ and
  $t'\in(\sigma_{\theta}(y'),s)$ such that $t'\in(t/2,t)$ since
  $\Omega_{\theta}\neq\emptyset $. It follows from \eqref{eq:tilde-sigma-def}
  that $|y'-y|>2Ls$ so $|x-y|\approx_{L} |y'-y|$ so for all $z\in\R$
  \begin{align*}
      &    sW_{s}(x-y)W_{t}(z-y)\lesssim_{L}sW_{s}(y'-y)W_{s}(z-y)\lesssim W_{t'}(z-y').
  \end{align*}

  Thus, since in both cases $(y',\xi+\theta t'^{-1},t')\in \Omega$ we have by the
  definition \eqref{eq:aux-map} of $\mb M$ that
  \begin{align*}
    sW_{s}(x-y) \mb M(y,\xi+\theta t^{-1},t) \lesssim_{L} \mb
    M(y',\xi+\theta t'^{-1},t)\leq \lambda
  \end{align*}
as required.
\end{proof}

The next technical lemma will be used as a continuous stopping time
argument. It relates the Lipschitz assumption on enlarged convex
regions of the previous statement with a crucial measurability estimate.

\begin{lem}[Continuous stopping time]\label{lem:continuous-stopping-time}
  Let  $\sigma:\R\to\R^{+}$ be a Lipschitz function with Lipschitz
  constant $L^{-1}<1$. Then the function 
  \begin{align*}
    \rho_{\sigma}(z):=\int_{\R}\frac{1}{2\sigma(x)}\1_{B_{\sigma(x)}}(z-x) \dd x
  \end{align*}
  satisfies
 $
    \left(1+\frac{2
        }{L-1}\right)^{-1}\!\!\!\!\!<\rho_{\sigma}(z)< 1+\frac{2 }{L-1} 
  $
  and in particular for any non-negative function $h(z)$ the 
  bounds 
  \begin{align*}
    \int_{\R}h(z)\dd z \approx_{L} \!\!\int_{\R}\;
    \fint_{B_{\sigma(x)}(x)} h(z) \dd z \;\dd x.
  \end{align*}
  hold.
\end{lem}

\begin{proof}
  Since $\sigma$ is $L^{-1}$-Lipschitz, for any $z\in\R$  we have that
  \begin{align*}
    B_{\left(1+L^{-1}\right)^{-1} \sigma (z)}(z)\subseteq    \{x \sthat z\in B_{\sigma (x)}(x)\}\subseteq B_{\left(1-L^{-1}\right)^{-1} \sigma (z)}(z).
  \end{align*}
  By the same reason on $\{x \sthat z\in B_{\sigma (x)}(x)\}$ we have
  that
  \begin{align*}
    \left(1+L^{-1}\right)^{-1} \sigma (z)\leq \sigma (x)\leq \left(1-L^{-1}\right)^{-1} \sigma (z).
  \end{align*}
  The conclusion follows.
\end{proof}


\subsection{Proof of Proposition \ref{prop:size-control}}

Let $T=T(x,\xi,s)$ be a tent and suppose that $K$ and $E$ are as in
\ref{prop:size-control}. Since the statement of Proposition
\ref{prop:size-control} is invariant under time and frequency
translations, we may assume, without loss of generality, that $T$ is
centered at the origin i.e. $T=T(0,0,s)$. If
$T\setminus (K\cup E)=\emptyset$  there is nothing to prove.
Let us set 
\begin{align*}
&  \Theta_{*}=\{\theta\in\Theta \sthat \exists(y,\theta t^{-1},t)\in  T(0,0,s) \setminus(K\cup E)\},
  \\
&\begin{aligned}
  &\Theta^{(i)}_{*}:=\Theta^{(i)}\cap \Theta_{*}&&\Theta^{(e)}_{*}:=\Theta^{(e)}\cap \Theta_{*}.
  \end{aligned}
\end{align*}
For $\theta\in\Theta_{*}$, using Lemma \ref{lem:convex-region-extention} we may assume that there
exists a $L^{-1}$-Lipschitz function $\sigma_{\theta}:\R\to(0,2s]$, with
$L>4$ sufficiently large to be chosen later, 
that satisfies condition \eqref{eq:tilde-sigma-condition-1} such that
\begin{align*}\num\label{eq:M-bound-in-proof}
&  s W_{s}(y)\mb M(y,\theta t^{-1},t) \lesssim \lambda &&\forall
                                                     y\in\R,\,\theta\in\Theta,\, t\in(\sigma(y),3s).
\end{align*}
Let us set  $\Omega=\bigcup_{\theta\in\Theta}\Omega_{\theta}$,
$\Omega^{(i)}=\bigcup_{\theta\in\Theta^{(i)}}\Omega_{\theta}$, and $\Omega^{(e)}=\bigcup_{\theta\in\Theta^{(e)}}\Omega_{\theta}$ with 
\begin{align*}
  \Omega_{\theta}=\left\{
  \begin{aligned}
&    \{(y,\theta t^{-1},t)\in T(0,0,s) \sthat
t>\sigma_{\theta}(y)\}&&\theta\in\Theta_{*}\\
&\emptyset&&\text{otherwise.}
    \end{aligned}
  \right.
\end{align*}
We need to show that
\begin{align*}
&  \|\mb A\, \1_{\mb X\setminus(K\cup E)}\|_{S_{m}(T)} \leq   \|\mb A\, \1_{\Omega}\|_{S_{m}(T)}\lesssim \lambda && \forall T\in\mb T
\end{align*}
or equivalently (see \eqref{eq:def:size-m}) that
\begin{align*}
  \|\mb A\,\1_{\Omega^{(i)}}\|_{S^{1}(T^{(i)})}\lesssim \lambda&&      \|\mb A\,\1_{\Omega}\|_{S^{2}(T)}\lesssim \lambda.
\end{align*}
In this proof all our implicit constants depend on the choice of $L$.

Let us fix a choice of left truncated wave packets
$\Psi_{y,\eta,t}^{c_{-},c_{+}}(z)$ in the defining expression
\eqref{eq:def:var-mass-embedding}. We will show that the statement
holds in this case. The proof for right truncated wave packets is
simmetric.

Comparing the definitions \eqref{eq:def:var-mass-embedding} and
\eqref{eq:aux-map} for $\mb A$ and $\mb M$ respectively, it follows
from the bound
$\left|\Psi_{y,\eta,t}^{\c_{k}(z),\c_{k+1}(z)}(z)\right|\leq
W_{t}(z-y)$ that
\begin{align*}\num\label{eq:S-lambda-infty}
  &  \mb A(y,\eta,t)\lesssim \mb M(y,\eta,t),
  &
  & \|\mb A \,\1_{\Omega}\|_{S^{\infty}(T)}:=
    \sup_{(y,\eta,t)\in
    \Omega}\mb A(y,\eta,t)\lesssim \lambda.
\end{align*}
This implies 
\begin{align*}\num\label{eq:Sm-L2-L1Linfty-dom}
  \frac{1}{s}\iiint_{\substack{(y,\eta,t)\in\Omega\\\eta<0}} |\mb
  A(y,\eta,t)|^{2} \lesssim \lambda  \;\frac{1}{s} \iiint_{\substack{(y,\eta,t)\in\Omega\\\eta<0}} |\mb
  A(y,\eta,t)|\dd y\dd \eta\dd  t
\end{align*}
and thus we may assume that
$\alpha^{-}=\beta^{-}<0<\beta^{+}<\alpha^{+}$ and we can reduce to showing
\begin{align*}\num\label{eq:S-lambda-bound}
  \|\mb A\,\1_{\Omega^{(i)}}\|_{S^{1}(T^{(i)})}\lesssim \lambda&&      \|\mb A\,\1_{\Omega^{(e)}}\|_{S^{2}(T^{(e)})}\lesssim \lambda.
\end{align*}

\subsubsection{Proof of the first inequality of \eqref{eq:S-lambda-bound}}
~

It holds that
\begin{align*}
&  \|\mb A\, \1_{\Omega^{(i)}}\|_{S^{1}(T^{(i)})}\approx \iiint_{\Omega^{(i)}} \left|\int_{\R}\sum_{k\in\Z}\a_{k}(z)\Psi_{y,\eta,t}^{\c_{k}(z),\c_{k+1}(z)}(z)
  \dd z \right| \dd y \dd \eta \dd t\\
&  \leq \frac{1}{s}\int_{\theta\in \Theta_{*}^{(i)}}
  \int_{y\in B_{s}} \int_{t=\sigma_{\theta}(y)}^{s} \int_{z\in\R}\sum_{k\in\Z} |\a_{k}(z)|
  \left|\Psi_{y,\theta t^{-1},t}^{\c_{k}(z),\c_{k+1}(z)}(z)\right|\dd z \frac{\dd
  t}{t}\dd y \dd \theta.
\end{align*}
According to \eqref{eq:eq:wave-packet-geomtery-vanish} The wave-packet
$\Psi_{y, \theta t^{-1}, t}^{\c_{k}(z),\c_{k+1}(z)}(z)$ vanishes unless
$\theta-t\c_{k}(z)\in B_{\epsilon}(d)$ and $t\c_{k+1}(z)
-\theta>d'$. Since $\theta\in\Theta^{(i)}\subset[-d',d-\epsilon]$, the
integrand vanishes unless $\c_{k}(z)<0<\c_{k+1}(z)$. Let
$k^{*}_{z}\in\Z$ be the index, if it exists, such that this inequality
holds and set $a^{*}(z):=\a_{k^{*}_{z}}(z)$, $c^{*}(z)=\c_{k^{*}(z)}(z)$. If no
such index exists simply set $a^{*}(z)=0$. 

Using that given $t<s$ and $y\in B_{s}$ one has 
\begin{align*}
  \Big|\Psi_{y,\theta t^{-1},t}^{\c_{k}(z),\c_{k+1}(z)}(z)\Big|\lesssim s W_{s}(z)\,
t\, W_{t}(z-y)^{2}\leq W_{t}(z-y)
\end{align*}
and using the statement of Lemma
\ref{lem:continuous-stopping-time} we have that 
\begin{align*}
  &\begin{aligned}\|\mb A \,\1_{\Omega^{(i)}}\|_{S^{1}(T^{(i)})}
    \lesssim 
    \frac{1}{s}
    \int_{\theta\in\Theta_{*}^{(i)}}\!\int_{y\in B_{s}}\!\int_{t=\sigma_{\theta}(y)}^{s}
    &\!\int_{x\in\R}
 \!\fint_{z \in B_{\sigma_{\theta}(x)}(x)} \hspace{-2em}|a^*(z)| sW_{s}(z)\,t\,W_{t}(y-z)^{2}
  \\
  &
    \hspace{-2.5em}\times
    \1_{B_{\epsilon}(d)}(\theta-tc^*(z)) \dd z \,\dd x \;\frac{\dd
    t}{t}\dd y \dd \theta
    = I+II
  \end{aligned}
  \\
  \shortintertext{where}
  &
    \begin{aligned}
      I:= \frac{1}{s}\int_{\theta\in \Theta_{*}^{(i)}} \int_{x\in\R}
      \fint_{z\in B_{\sigma_{\theta}(x)}(x)} \hspace{-2.5em}|a^*(z)|
      \int_{t=(1-2/L)\sigma_{\theta}(x)}^{s} &\;\int_{y\in
        B_{s}}W_{t}(y-z) \dd y \\&\hspace{-1em}\times
      \1_{B_{\epsilon}(d)}(\theta-tc^*(z)) \frac{\dd t}{t} \;\dd z
      \dd x \dd \theta
    \end{aligned}\\
  &
    \begin{aligned}
      II:= \frac{1}{s}\int_{\theta\in \Theta_{*}^{(i)}}
      \int_{x\in\R} \int_{y\in B_{s}}
      \int_{t=\sigma_{\theta}(y)}^{(1-2/L)\sigma_{\theta}(x)}
      \!\!&\fint_{z\in B_{\sigma_{\theta}(x)}(x)} \hspace{-2.5em}|a^*(z)|\,s\,W_{s}(z)\,t\,
      W_{t}(y-z)^{2}\\&\hspace{2.5em}\times \1_{B_{\epsilon}(d)}(\theta-tc^*(z))
      \dd z \;\frac{\dd t}{t} \dd y \dd x \dd \theta
    \end{aligned}
\end{align*}
Suppose that $L>2\frac{\alpha^{+}-\alpha^{-}}{\alpha^{+}-d-\epsilon}$ so
that for any $c\in\R$ one has
\begin{align*}\num\label{eq:L-choice-condition}
  &\left\{
    \begin{aligned}
      &\theta-t c\in B_{\epsilon}(d)    \\
      &t>(1-2/L)\sigma_{\theta}(x)  \end{aligned}
\right. &&\implies && \theta-\sigma_{\theta}(x) c\in\Theta.
\end{align*}
We begin by estimating the term $I$. Notice that if
$|x|>2Ls$ then integrand vanishes. We bound $I$ by the
auxiliary embedding map \eqref{eq:aux-map} as follows:
\begin{align*}
  I  \lesssim&
     \frac{1}{s}\int_{\theta\in \Theta_{*}^{(i)}}
    \int_{x\in \R}  
   \fint_{z \in B_{\sigma_{\theta}(x)}(x)}
    \hspace{-2em}|a^*(z)|  \int_{t=(1-2/L)\sigma_{\theta}(x)}^{s} \hspace{-4em}\1_{B_{\epsilon}(d)}(\theta-tc^*(z))
    \frac{\dd t}{t}\dd z \dd x \dd \theta
  \\\lesssim&
     \frac{1}{s}\int_{\theta\in \Theta_{*}^{(i)}}
    \int_{x\in B_{2Ls}} \fint_{B_{\sigma_{\theta}(x)}(x)}
    \hspace{-2em}|a^*(z)|\, \1_{\Theta_{*}}
      \left(\theta-\sigma_{\theta}(x)\c_{k}(z)\right)
              \ln\left(\frac{\beta^{+}-\theta+3\epsilon}{\beta^{+}-\theta}\right)
              \,\dd z\dd x\,\dd \theta
      \\\leq&
   \frac{1}{s}   \int_{\theta\in \Theta_{*}^{(i)}} 
      \int_{x\in B_{2Ls}} \mb M
          (x,\theta\sigma_{\theta}(x)^{-1},\sigma_{\theta}(x)) \dd x
         \, \ln\left(\frac{\beta^{+}-\theta+3\epsilon}{\beta^{+}-\theta}\right)
          \dd \theta
      \lesssim \lambda.
\end{align*}
The last inequality holds since
$(x,\theta\sigma_{\theta}(x)^{-1},\sigma_{\theta}(x))\in\Omega_{\theta}$
and follows from \eqref{eq:M-bound-in-proof}.

We now estimate the term $II$. Notice that
\begin{align*}\num\label{eq:off-scale-condition}
  (1-2/L)\sigma_{\theta}(x)>\sigma_{\theta}(y) \implies
  |x-y|\geq L(\sigma_{\theta}(x)-\sigma_{\theta}(y))> 2\sigma_{\theta}(x)
\end{align*}
Thus if 
$z \in B_{\sigma_{\theta}(x)}(x)$ then $|y-z|>\sigma_{\theta}(x)>t$, 
$|x-y|\approx |y-z|$, and also $s W_{s}(z)\lesssim sW_{s}(x)$ so 
\begin{align*}
&\begin{aligned}  II\lesssim  \int_{\theta\in \Theta_{*}^{(i)}} \int_{x\in \R}\hspace{-1em}W_{s}(x) 
      \int_{y\in B_{s}}     
      &\int_{t=\sigma_{\theta}(y)}^{(1-2/L)\sigma_{\theta}(x)} t \,W_{t}(y-x)
  \\
      \times& \fint_{B_{\sigma_{\theta}(x)}(x)}\hspace{-2em} |a^*(z)|\,W_{t}(y-z)
      \,\1_{B_{\epsilon}(d)}(\theta-tc^{*}(z)) \dd z \frac{\dd t}{t}
      \dd y \dd x \dd \theta
    \end{aligned}
  \\&\phantom{II} \lesssim
      \int_{\theta\in\Theta_{*}^{(i)}}
      \int_{ x\in \R}   \hspace{-1em}W_{s}(x)   \int_{y\in B_{s}}
              \int_{t=\sigma_{\theta}(y)}^{(1-2/L)\sigma_{\theta}(x)} \!\!\frac{t}{2\sigma_{\theta}(x)}
      W_{t}(y-x)  \mb M(y,\theta t^{-1},t)\frac{\dd t}{t}\dd y \dd x 
              \dd \theta.
\end{align*}
Since the inmost integral vanishes unless
$|y-x|>2\sigma_{\theta}(x)$,  we have that
\begin{align*}
   \int_{t=\sigma_{\theta}(y)}^{(1-2/L)\sigma_{\theta}(x)}\!\! \frac{t}{2\sigma_{\theta}(x)}
      W_{t}(y-x) \frac{\dd t}{t} \lesssim W_{\sigma_{\theta}(x)}(y-x)
\end{align*}
and so using \eqref{eq:M-bound-in-proof} we obtain
\begin{align*}II\lesssim&
          \lambda\int_{\theta\in\Theta_{*}^{(i)}}\!
      \int_{x\in \R}W_{s}(x) \!\int_{y\in B_{s}}
          W_{\sigma_{\theta}(x)}(y-x) \dd y \dd x
                          \,\dd \theta
  \\\lesssim
&\lambda\int_{\theta\in\Theta_{*}^{(i)}}\!
      \int_{x\in \R}  W_{s}(x)  \,\dd x\,
          \dd \theta \lesssim\lambda.
\end{align*}
This concludes the proof for the first bound of
\eqref{eq:S-lambda-bound}.


\subsubsection{Proof of the second inequality of \eqref{eq:S-lambda-bound}}

As noted in \eqref{eq:Sm-L2-L1Linfty-dom} we may suppose that
$\Theta^{(e)}=[\beta^{+},\alpha^{+})$; for ease of notation set
$\Omega^{(e)}=\Omega\cap T^{(e)}$ so the required quantity to bound becomes
$\|\mb A\,\1_{\Omega^{(e)}}\|_{S^{2}(T^{(e)})}$.  We concentrate on showing
the dual bound
\begin{align*}\num\label{eq:dim:sm2-dualizing}
  &\frac{1}{s}\Bigg|\iiint_{\mb X}
    h(y,\eta,t) \int_{\R} \sum_{k\in\Z} \a_{k}(z)
    \Psi_{y,\eta,t}^{\c_{k}(z),\c_{k+1}(z)}(z)\dd z\,  \dd y\dd \eta
    \dd t \Bigg|\lesssim \lambda \frac{\|h(y,\eta,t)\|_{L^{2}}}{s^{1/2}}
\end{align*}
for any $h\in C^{\infty}_{c}(\Omega^{(e)})$ where the $\|\cdot\|_{L^{2}}$
is the classical Lebesgue $L^{2}$ norm relative to the measure $\dd
y\dd \eta\dd t$.  A change of variables and
the Minkowski
inequality give 
\begin{align*}
&  \frac{1}{s} \Bigg|\iiint_{\mb X}
  h(y,\eta,t) \int_{\R} \sum_{k\in\Z} \a_{k}(z)
    \Psi_{y,\eta,t}^{\c_{k}(z),\c_{k+1}(z)}(z)\dd z\,  \dd y\dd \eta
  \dd t \Bigg|\\
  &\leq 
  \frac{1}{s} \int_{\R}\sum_{k\in\Z}|\a_{k}(z)|
    \left|\iiint_{\Omega^{(e)}} h(y,\eta,t)
    \Psi_{y,\eta,t}^{\c_{k}(z),\c_{k+1}(z)}(z) \dd  y  \dd
      \eta \dd t\right| \dd z   
      \\
&\leq \frac{1}{s}\int_{\theta\in\Theta_{*}^{(e)}}  \int_{\R}\sum_{k\in\Z}|\a_{k}(z)|
      \left|\int_{y\in B_{s}} \!\int_{t=\sigma_{\theta}(y)}^{s}  \hspace{-2em}h(y,t^{-1}\theta,t)
        \Psi_{y,\theta t^{-1},t}^{\c_{k}(z),\c_{k+1}(z)}(z) \dd  y \frac{\dd t}{t} \right| \dd z   \dd
     \theta.
\end{align*}
On the other hand the Hölder inequality gives that
  \begin{align*}
    &\int_{\theta\in\Theta_{*}^{(e)}}\|h(y,t^{-1}\theta,t)\|_{L^{2}_{(\dd y
      \dd t/t)}} \dd \theta
      = \int_{\theta\in\Theta_{*}^{(e)}}\!\! \Big(\int_{y\in B_{s}} \!\int_{t=\sigma_{\theta}(y)}^{s}  \hspace{-2em} |h(y,t^{-1}\theta,t)|^{2} \dd y \frac{\dd
    t}{t} \Big)^{\frac{1}{2}} \dd \theta\\
    & \lesssim
    \|h(y,\eta,t)\|_{L^{2}}
  \end{align*}
  where $\|\cdot\|_{L^{2}_{(\dd y \dd t/t)}}$ is the classic
  Lebesgue $L^{2}$ norm with respect to the measure $\dd y\frac{\dd t}{t}$.
  Thus \eqref{eq:dim:sm2-dualizing} follows by showing 
  \begin{align*}\num\label{eq:dim:sm2-dual-estimate}
    \frac{1}{s} \int_{\R}\sum_{k\in\Z}|\a_{k}(z)|
    \Bigg|\int_{y\in B_{s}} \!\int_{t=\sigma_{\theta}(y)}^{s}
    \hspace{-2em} h(y,t^{-1}\theta,t)
    &\Psi_{y,t^{-1}\theta,t}^{\c_{k}(z),\c_{k+1}(z)}(z) \dd  y \frac{\dd
      t}{t} \Bigg| \dd z
    \\&
        \lesssim \lambda
        \frac{\left\|h(y,t^{-1}\theta,t)\right\|_{L^{2}_{(\dd  y\dd t/t)}}}{s^{1/2}}
  \end{align*}
  with a constant uniform in $\theta\in\Theta_{*}^{(e)}$.
  For sake of notation from now on we will omit  the dependence on $\theta$ by
  writing
  \begin{align*}
    &h(y,t):=h(y,t^{-1}\theta,t)
    &&\Psi_{y,t}^{c_{-},c_{+}}(z):=\Psi_{y,t^{-1}\theta,t}^{c_{-},c_{+}}(z)&& \sigma(x):=\sigma_{\theta}(x). 
  \end{align*}
  Using the above notation and Lemma

  \ref{lem:continuous-stopping-time} we write
\begin{align*}
  &\frac{1}{s} \int_{\R} \sum_{k\in\Z}|\a_{k}(z)|
    \Bigg|\int_{y\in B_{s}} \!\int_{t=\sigma(y)}^{s}
    \hspace{-2em} h(y,t)
    \Psi_{y,t}^{\c_{k}(z),\c_{k+1}(z)}(z) \dd  y \frac{\dd
    t}{t} \Bigg| \dd z \lesssim I+II
  \\
  \shortintertext{where}
  I:&= \frac{1}{s} \int_{x\in B_{2Ls}}\fint_{z\in B_{\sigma(x)}(x)}\sum_{k\in\Z}|\a_{k}(z)|
    \Bigg|\int_{y\in B_{s}} \!\int\mcl_{\qquad t=(1-2/L)\sigma(x)}^{s} h(y,t)
    \Psi_{y,t}^{\c_{k}(z),\c_{k+1}(z)}(z) \dd  y \frac{\dd
    t}{t} \Bigg| \dd z  \dd x
  \\
  II:&= \frac{1}{s} \int_{x\in \R}\fint_{z\in B_{\sigma(x)}(x)}\sum_{k\in\Z}|\a_{k}(z)|
    \Bigg|  \int_{y\in B_{s}} \!\int\mcl_{\qquad t=\sigma(y)}^{\mathclap{\quad(1-2/L)\sigma(x)}}  h(y,t)
    \Psi_{y,t}^{\c_{k}(z),\c_{k+1}(z)}(z) \dd  y \frac{\dd
    t}{t} \Bigg| \dd z   \dd x
\end{align*}

We start with bounding $I$. Suppose that $L>1$ is chosen large
enough so that \eqref{eq:L-choice-condition} holds and recall that
$\Psi_{y,t}^{\c_{k}(z),\c_{k+1}(z)}(z)=0$ unless $\theta-t\c_{k}(z)\in
B_{\epsilon}(d)$. We thus have 
\begin{align*}
  I  =&\frac{1}{s} \int_{x\in B_{2Ls}}\fint_{z\in B_{\sigma(x)}(x)}\sum_{k\in\Z}|\a_{k}(z)|\,
        \bigg|\!\int_{y\in B_{s}} \!\int\mcl_{\qquad t=(1-\frac{2}{L})\sigma(x)}^{s} h(y,t)
        \Psi_{y,t}^{\c_{k}(z),\c_{k+1}(z)}(z) \dd  y \frac{\dd
        t}{t} \bigg| \dd z  \dd x
  \\\leq& \frac{1}{s}\int_{x\in B_{2Ls}}
          \fint_{B_{\sigma(x)}(x)} \!\!\Big(\!
          \sum_{k\in\Z}
          |\a_{k}(z)|^{r'}\1_{\Theta}\big(\theta-\sigma(x)
          \c_{k}(z)\big) \Big)^{1/r'} 
          \mathcal{H}_{x}(z)
          \dd z\, \dd x\\
  \leq& \frac{1}{s}\int_{x\in B_{2Ls}}
        \mb M(x,\theta\sigma(x)^{-1},\sigma(x))
        \sup_{z\in        B_{\sigma(x)}(x)} \mc H_{x}(z) \dd x\leq \frac{\lambda}{s}\int_{x\in B_{2Ls}}
        \sup_{\mathclap{\qquad \;z\in B_{\sigma(x)}(x)}}\; \mc H_{x}(z) \dd x,\\
  \shortintertext{where}
      &  \mc H_{x}(z):=\bigg(\sum_{k\in\Z}\bigg|\!\int_{y\in B_{s}} \!\int\mcl_{\qquad t=(1-\frac{2}{L})\sigma(x)}^{s} h(y,t)
        \Psi_{y,t}^{\c_{k}(z),\c_{k+1}(z)}(z) \dd  y \frac{\dd
        t}{t} \bigg|^{r}\bigg)^{1/r}.
\end{align*}
We claim that 
\begin{align*}\num\label{eq:Hxz-domination}
  &  \mc H_{x}(z)\lesssim \mc V_{\sigma(x)}^{r} H_{s} (z) + \mc E_{\sigma(x)}(z)    &&H_{\tau}(z):=\int_{t=0}^{\tau}\int_{y\in B_{s}} h(y,t) \Psi_{y,t}^{0,+\infty}(z) \dd
    y\frac{\dd t}{t}\\
  &\mc E_{\sigma(x)}(z) := \left(\int_{(1-2/L)\sigma(x)}^{s}
    \hspace{-2.5em}|h^{*}(z,t)|^{2}\frac{\dd t}{t}\right)^{1/2}
                                                                       &&h^{*}(z,t):=\int_{\R}|h(y,t)|W_{t}(z-y)\dd y
\end{align*}
with $\mc V^{r}_{\sigma(x)}$ defined in Lemma
\ref{lem:var-trunc}, and that
\begin{align*}
\num\label{eq:Hxz-domination-L2-E}
  &\|\mc E_{0}\|_{L^{2}}\lesssim
  \|h\|_{L^{2}_{(\dd y\dd t/t)}}&
                               & \sup_{z\in B_{\sigma(x)}(x)} \mc
                                 E_{\sigma(x)}(z)\lesssim \fint_{B_{2\sigma(x)}(z)}    \mc
                                 E_{0}(z) \dd z
  \\\num\label{eq:Hxz-domination-L2-H}
  &\|H_{s}\|_{L^{2}}\lesssim \|h\|_{L^{2}_{(\dd y\dd t/t)}}&&\sup_{z\in B_{\sigma(x)}(x)} \mc H_{x}(z) \dd x \lesssim
    \fint_{B_{2\sigma(x)}(x)}M \mc V^{r} H_{s}(z) \dd z.
\end{align*}
This would provide us with the required bounds for $I$. As a matter of
fact, according to Lemma \ref{lem:continuous-stopping-time} and
\ref{lem:var-trunc} we  have that
\begin{align*}
  &
    I\lesssim\!
    \frac{\lambda}{s}  \int\limits_  {x\in
    B_{2Ls}}  \!\!\fint\limits_ {\;\;B_{2\sigma(x)}(x)}\!\!\!\Big(\!M \mc V^{r} H_{s}(z) + \mc
                                 E_{\sigma(x)}(z) \Big) \dd z\lesssim_{R}  \!  \frac{\lambda}{s^{1/2}}\Big(
    \|M \mc V^{r} H_{s}\|_{L^{2}}+\|\mc E_{0}(z)\|_{L^{2}}\Big)
  \\
  &\lesssim  \frac{\lambda}{s^{1/2}}\Big(\|H_{s}(z)\|_{L^{2}}+\|\mc
    E_{0}(z)\|_{L^{2}}\Big)\lesssim \lambda \frac{\|h\|_{L^{2}_{(\dd y\dd t/t)}}}{s^{1/2}}
\end{align*}
as required.

The first bound of \eqref{eq:Hxz-domination-L2-E} follows by the Young
inequality and Fubini:
\begin{align*}
&  \|\mc E_{0}\|_{L^{2}} \leq \iint_{\R\times\R^{+}}\Big|\int_{\R} |h(y,t)|W_{t}(z-y)\dd
  y\Big|^{2}\dd z\frac{\dd t}{t}\\&\leq
                                     \int_{\R^{+}}\int_{\R} |h(y,t)|^{2}\dd y\; \Big(\int
  W_{t}(z-y)\dd y\Big)^{2} \frac{\dd t}{t} \lesssim \|h\|_{L^{2}_{(\dd
  y\dd t/t)}}^{2}.
\end{align*}
The second bound follows from the fact that for
small enough $\epsilon>0$ and as long as $|z-z'|<\epsilon t$ the bound
\begin{align*}
&  |h^{*}(z,t)- h^{*}(z',t)| \leq \int_{\R} |h(y,t)|
  |W_{t}(z-y)-W_{t}(z'-y)|\dd y
  \\&\leq  2^{-100}\int_{\R} |h(y,t)|
   W_{t}(z-y)\dd y=2^{-100} h^{*}(z,t)
\end{align*}
holds so similarly
\begin{align*}
  \Big|\mc E_{\sigma(x)}(z) -\mc E_{\sigma(x)}(z') \Big|\leq 2^{-100} \mc E_{\sigma(x)}(z)
\end{align*}
as long as $|z-z'|<\epsilon \sigma(x)$ for some sufficiently small $\epsilon>0$.
\begin{align*}
  \fint_{B_{\sigma(x)}(z)}\mc E_{0}(z')\dd z' \gtrsim
  \fint_{B_{\sigma(x)}(z)}\mc E_{\sigma(x)}(z')\dd z' \gtrsim_{\epsilon}
  \fint_{B_{\epsilon\sigma(x)}(z)}\mc E_{\sigma(x)}(z)\dd z' = \mc E_{\sigma(x)}(z)
\end{align*}
and  the claim follows.

The first bound of \eqref{eq:Hxz-domination-L2-H} uses standard oscillatory
integral techniques: notice that for $t>t'$ one has
\begin{align*}
  \Big|\int_{\R} \Psi_{y,t}(z)  \overline{\Psi_{y',t'}(z)}\dd z\Big|
  \lesssim \frac{t'}{t} W_{t}(y-y') 
\end{align*}
so
\begin{align*}
\int|H_{s}(z)|^{2}(z) \lesssim 2 \int_{t=0}^{s}\int_{t'=0}^{t} \int_{y\in
  B_{s}}\int_{y'\in B_{s}}\hspace{-1em}| h(y,t)||h(y',t')| W_{t}(y-y') \dd y\dd
  y' \frac{t'}{t}\frac{\dd t}{t}\frac{\dd t'}{t'}&\\
\lesssim \|h\|_{L^{2}_{(\dd y\dd t/t)}}^{2}.&
\end{align*}
The second bound follows directly from Lemma \ref{lem:var-trunc}.

It remains to show inequality \eqref{eq:Hxz-domination}.
Notice that 
\begin{align*}
\mathcal{H}_{x}(z)=      \Bigg(\sum_{\substack{k\in\Z}}
  \Big|\int_{t_{k}^{-}(z)}^{t_{k}^{+}(z)} \int_{y\in B_{s}} h(y,t)
    \Psi_{y,t}^{\c_{k}(z),\c_{k+1}(z)}(z) \dd y
  \frac{\dd t}{t}\Big|^{r}\Bigg)^{1/r}
\end{align*}
where for  $k\in\Z$ we set 
\begin{align*}\num\label{eq:t+-}
 &  t_{k}^{+}(z) :=\sup\left\{t\in\big((1-2/L)\sigma(x),s\big) \sthat
    \Psi_{y,t}^{\c_{k}(z),\c_{k+1}(z)}(z) \neq 0\right\} 
  \\
  & t_{k}^{-}(z) := \inf\left\{t\in\big((1-2/L)\sigma(x),s\big) \sthat
    \Psi_{y,t}^{\c_{k}(z),\c_{k+1}(z)}(z) \neq 0\right\}.
\end{align*}
We have omitted writing the implicit dependence on $x\in\R$ and we
will simply ignore the indexes $k\in\Z$ for which the above sets
are empty. Notice that the intervals $\big[t_{k}^{-}(z),\,t_{k}^{+}(z)\big)$
are disjoint.  According to the conditions
\eqref{eq:eq:wave-packet-geomtery-vanish} on the geometry of truncated wave
packets the following bounds hold:
\begin{align*}\num\label{eq:t-bounds}
  &t_{k}^{+}(z)\c_{k}(z)\in \overline{B_{\epsilon}(\theta-d)}
  &&t_{k}^{-}(z)\c_{k+1}(z)\geq \theta+d' .
\end{align*}
Using the  smoothness conditions 
\eqref{eq:wave-packet-smoothness} on the wave packets and writing a
Lagrange remainder term we have that
\begin{align*}\num\label{eq:wave-packet-remainder}
  \left|\Psi_{y,t}^{\c_{k}(z),\c_{k+1}(z)}(z) -
  \Psi_{y,t}^{0,+\infty}(z)\right|\leq  \Big(|t\c_{k}(z)|+
  \max\left(d''-\theta-t\c_{k+1}(z);0\right) \Big) W_{t}(y-z)
\end{align*}
so the bound 
\begin{align*}
  \mc H_{x}(z) &\leq \mc H_{x,1}(z) + \mc H_{x,2}(z)
  \\ \mc H_{x,1}(z) &:=\Big(\sum_{k\in\Z} \Big|\int_{y\in
                     B_{s}}\int_{ t=t_{k}^{-}(z)}^{t_{k}^{+}(z)} h(y,t)
                     \Psi_{y,t}^{0,+\infty}(z) \dd y
                     \frac{\dd t}{t}\Big|^{r}\Big)^{1/r}
  \\
              &=\Big(\sum_{k\in\Z} |H_{t_{k}^{+}(z)}(z)-H_{t_{k}^{-}(z)}(z)|^{r}\Big)^{1/r}
  \\
  \mc H_{x,2}(z) &:=
                  \Big(\sum_{k\in\Z} \Big|\int_{t_{k}^{-}(z)}^{t_{k}^{+}(z)} h^{*}(z,t)
  \,\big(|t\c_{k}(z)|+  \max(d''-\theta-t\c_{k+1}(z),0)\big)
  \frac{\dd t}{t}\Big|^{r}\Big)^{1/r}
\end{align*}
holds.  Notice that
\begin{align*}
&  \int_{t_{k}^{-}(z)}^{t_{k}^{+}(z)} t^{2}|\c_{k}(z)|^{2}\frac{\dd t}{t}\leq
  \frac{|t^{+}_{k}(z)\c_{k}(z)|^{2}}{2}\leq C_{\alpha^{+}}
  \\
&
  \int_{t_{k}^{-}(z)}^{t_{k}^{+}(z)} \max\big(d''-\theta-t\c_{k+1}(z);0\big) \frac{\dd t}{t}\leq
   \int_{t_{k}^{-}\c_{k+1}(z)}^{d''-\theta}\frac{\dd t}{t} \leq
  \int_{d'+\theta}^{d''-\theta}\frac{\dd
  t}{t}\lesssim C_{d',d,\alpha^{+},\beta+}
\end{align*}
for some constant $C_{\alpha^{+}}$ and $C_{d',d,\alpha^{+},\beta+}$.
Since $r>2$, Cauchy-Schwartz gives
\begin{align*}
  \mc H_{x,2}(z) \leq \left(
  \sum_{k\in\Z}\int_{t_{k}^{-}(z)}^{t_{k}^{+}(z)}
  h^{*}(z,t)^{2} \frac{\dd t} {t} \right)^{\frac{1}{2}}\lesssim  \mc E_{\sigma(x)}(z).
\end{align*}
This is consistent with \eqref{eq:Hxz-domination}.
To estimate $\mc H_{x,1}$: introduce a frequency cutoff $\Upsilon\in
S(\R)$ such that
\begin{align*}
  &\FT \Upsilon\in C^{\infty}_{c}(B_{\theta+b})
  &&\FT \Upsilon\geq 0
  &&\FT \Upsilon= 1\text{ on } B_{\theta}&&
  \Upsilon_{\tau}(z):=\tau^{-1}\Upsilon\left(\frac{z}{\tau}\right).
\end{align*}
According to 
\eqref{eq:wave-packet-freq-support},  $\FT \Psi_{y,t}^{0,+\infty}$ is supported on
$B_{t^{-1}b}(t^{-1}\theta)$ so one has the following
  \begin{align*}
      &\Psi_{y,t}^{0,+\infty}*\Upsilon_{\tau}(z) =\Psi_{y,t}^{0,+\infty}(z)
      &&   \text{if }   \frac{t}{\tau}\geq\frac{\theta+b}{\theta}\\
      &\Psi_{y,t}^{0,+\infty}*\Upsilon_{\tau}(z) =0
      &&\text{if } \frac{t}{\tau}<\frac{\theta-b}{\theta+b}\\
      &\big|\Psi_{y,t}^{0,+\infty}*\Upsilon_{\tau}(z)\big|\lesssim
      W_{t}(z-y) && \text{if }
      \frac{\theta-b}{\theta+b}\leq\frac{t}{\tau}<\frac{\theta+b}{\theta}.
  \end{align*}
  Thus
  \begin{align*}
    &\Big|H_{\tau}-
      H_{s}*\Upsilon_{\tau}(z) \Big|\lesssim
      \int_{\tau\frac{\theta-b}{\theta+b}}^{\tau\frac{\theta+b}{\theta}}h^{*}(z,t)\frac{\dd t}{t}
  \end{align*}
  so
  \begin{align*}
    &\mc H_{x,1}(z) \lesssim 
      \Big(\sum_{k\in\Z}|H_{s}*\Upsilon_{t_{k}^{+}(z)}-H_{s}*\Upsilon_{t_{k}^{-}(z)}|^{r}\Big)^{\frac{1}{r}}
      +   \Big(\sum_{k\in\Z}\Big|\int_{t_{k}^{-}(z)\frac{\theta-b}{\theta+b}}^{t_{k}^{-}(z)\frac{\theta+b}{\theta}}h^{*}(z,t)\frac{\dd
      t}{t}\Big|^{r} \Big)^{\frac{1}{r}}
    \\
    &+
      \Big(\sum_{k\in\Z}\Big|\int_{t_{k}^{+}(z)\frac{\theta-b}{\theta+b}}^{t_{k}^{+}(z)\frac{\theta+b}{\theta}}h^{*}(z,t)\frac{\dd t}{t}\Big|^{r}
      \Big)^{\frac{1}{r}}\lesssim
      \mc V^{r}_{\sigma(x)}H_{s}(z) + \mc E_{\sigma(x)}(z)
  \end{align*}
  thus concluding the proof of \eqref{eq:Hxz-domination} and the bound
  on the term $I$.

  The estimate for the term $II$ can be done in a manner similar to
  the term $II$ in for the $S^{1}$ part of the size. Recall
  \eqref{eq:off-scale-condition}
  so that in expression for $II$ one has that 
  $z\in B_{\sigma(x)}(x)$, 
  $|y-z|>\sigma(x)>t$, and $|x-y|\approx|y-z|$. We also have that
  $y\in B_{s}$ so
  \begin{align*}
    \Psi^{\c_{k}(z),\c_{k+1}(z)}_{y,t}(z)\lesssim s W_{s} (x)\, t W_{t}(z-y)^{2}
  \end{align*}
  and $\Psi^{\c_{k}(z),\c_{k+1}(z)}_{y,t}(z)=0$ unless
  $t\c_{k}(z)<\theta<t\c_{k+1}(z)$, thus
  \begin{align*}
    II\lesssim&
                \frac{1}{s}\int_{x\in \R}sW_{s}(x)\int_{y\in
                B_{s}}\int_{t=\sigma (y)}^{(1-2/L)\sigma (x)}
                \hspace{-2em}tW_{t}(y-x) h(y,t)    
                \fint_{B_{\sigma (x)}(x)} \sum_{k\in\Z}|\a_{k}(z)|
    \\&\hspace{7em}\times
        \1_{\Theta}(\theta-t\c_{k}(z))\, W_{t}(z-y)
        \1_{(\c_{k}(z),\c_{k+1}(z))}(t^{-1}\theta) \dd z \frac{\dd
        t}{t}\dd y \dd x
    \\
    \lesssim& \int_{x\in \R}W_{s}(x) \int_{y\in B_{s}}\int_{t=\sigma
              (y)}^{(1x-2/L)\sigma
              (x)} \frac{t}{2\sigma(x)}W_{t}(y-x) h(y,t)
    \\&\times\int_{B_{\sigma (x)}(x)} \Big(\sum_{k\in\Z}|\a_{k}(z)|^{r'}
              \1_{\Theta}(\theta-t\c_{k}(z))\Big)^{1/r'}  W_{t}(y-z)\,\dd
        z\;\frac{\dd t}{t} \dd y \dd x\\
    \lesssim& \int_{x\in \R}W_{s}(x) \int_{y\in B_{s}}\int_{t=\sigma
              (y)}^{(1x-2/L)\sigma
              (x)} \frac{t}{2\sigma(x)}W_{t}(y-x) h(y,t) \mb M(y,\theta
              t^{-1},t)\frac{\dd t}{t} \dd y \dd x
  \end{align*}
Since the inmost integral vanishes unless
$|y-x|>2\sigma_{\theta}(x)$,  we have that
\begin{align*}
\Big(   \int_{t=\sigma_{\theta}(y)}^{(1-2/L)\sigma_{\theta}(x)}\big| \frac{t}{2\sigma(x)}
      W_{t}(y-x)\big|^{2} \frac{\dd t}{t} \Big)^{1/2}\lesssim W_{\sigma(x)}(y-x)
\end{align*}
it follows that 
  \begin{align*}
    II\lesssim& \lambda \int_{x\in \R}W_{s}(x) \int_{y\in B_{s}}W_{\sigma(x)}(y-x)
                \Big(\int_{t=0}^{s} |h(y,t)|^{2}\frac{\dd t}{t}\Big)^{1/2} \dd y \dd x
    \\
    \lesssim& \lambda\int_{x\in \R}W_{s}(x) M\bigg(\Big(\int_{t=0}^{s}
              |h(\cdot,t)|^{2}\frac{\dd t}{t}\Big)^{1/2}\bigg) (x) \dd x
    \\    
    \lesssim& \frac{\lambda\|h\|_{L^{2}_{(\dd y \dd t/t)}}} {s^{1/2}}
  \end{align*}
This concludes the proof.
  
  \qed

\section{The energy embedding and non-iterated
  bounds}\label{sec:e-bounds}
\subsection{The energy embedding}

Here we comment on how to deduce Theorem \ref{thm:energy-boundedness}
from the result in \cite{di2015modulation}. Let us fix a
$p\in (1,\infty]$ and $q\in\big(\max(2;p'),+\infty\big]$ and without loss of
generality let us suppose that $\FT f\in C^{\infty}_{c}(\R)$.
We will show that the weak versions of \eqref{eq:e-bounds} holds i.e.
\begin{align*}\num\label{eq:e-bounds-w}
  \|F\|_{L^{p,\infty}\L^{q}(S_{e})}\lesssim \|f\|_{L^{p}}.
\end{align*}
By interpolation this would allow us to conclude the strong bounds of
\eqref{eq:e-bounds}.

The paper \cite{di2015modulation} deals with embeddings into the space
$\mb X$ that they denote by $\mc Z$. The generating collection of
tents that they make use of is described in Section 2.1.2 of that
paper. Notice that the set of geometric parameters for the tents in
the present paper (Section \ref{sec:outer-measure-def}) is larger
than the one in \cite{di2015modulation} but a careful perusal of the
proofs therein shows that the same statements hold for the extended
range of parameters.

Let us recall the main statements from \cite{di2015modulation}.

\begin{thm}[Theorem 1 of \cite{di2015modulation}]\label{thm:thm-diou}
  Let $f\in S(\R)$ with $\FT f\in C^{\infty}_{c}$. Let $p\in(1, 2)$  and consider the set $\mc I_{f,\lambda,p}$ of
  maximal dyadic intervals contained in
  \begin{align*}\num\label{eq:maximal-superlevel}
    &\mc K_{f,\lambda,p}=\{x\in \R \sthat M_{p}f(x)>\lambda\}
    &
    &\text{and
      let}&&
      K_{f,\lambda ,p}:=\bigcup_{B_{\tau}(\zeta)\in \mc I_{f,\lambda,p}} D(\zeta,3\tau).
  \end{align*}
  Then with $q\in(p',\infty]$.
  \begin{align*}
    &\|F\1_{\mb X\setminus K_{f,\lambda,p}}\|_{L^{q}(S_{e})}\lesssim_{q,p}
  \lambda^{1-p/q} \|f\|_{L^{p}}^{p/q}.
\end{align*}
\end{thm}
We used the super level set of $M_{p}f$ instead of the super level set
of $M_{p}\big(Mf\big)$ to define $\mc K_{f,\lambda,p}$.  As mentioned
in section 7.3.1 of \cite{di2015modulation}, the inner maximal
function appears only in the reduction from the case with $\FT f$
compactly supported to the case with a general $f\in S(\R)$. By our
assumptions we can effectively ignore this complication.

\begin{prop}[Proposition 3.2 + equations (2.6) and (2.7) of \cite{di2015modulation}]\label{prop:local-L2}
  The estimate
  \begin{align*}
    &\|F\,\1_{D(x,s)} \|_{L^{q}(S_{e})} \lesssim_{N,q}
      \left(1+\frac{\dist(\spt f ; B_{s}(x)}{s}\right)^{-N}\|f\|_{L^{q}}
  \end{align*}
  holds for all  $N>0$ and $q\in(2,\infty]$.
\end{prop}

\begin{lem}[Equation (7.3) of
  \cite{di2015modulation}]\label{lem:locality-decay}
  The estimate
  \begin{align*}
    \|F\,\1_{D(x,s)} \|_{L^{\infty}(S_{e})} \lesssim_{N}
    \left(1+\frac{\dist\big(\spt f ; B_{s}(x)\big)}{s}\right)^{-N}
    \inf_{z\in B_{s}(x)} Mf(z)
  \end{align*}
  holds for any $N>0$.
\end{lem}

\begin{cor}\label{cor:local}
  Suppose that $\spt f\cap B_{2s}(x)=\emptyset$ then
  \begin{align*}
    \|F \1_{D(x,s)}\|_{L^{q}(S_{e})} s^{-1/q}\lesssim_{N,p}
    \left(1+\frac{\dist\big (\spt f ; B_{s}(x)\big)}{s}\right)^{-N} s^{-1/p}\|f\|_{L^{p}}
  \end{align*}
  for all $p\in[1,2)$, $q>p'$, and $N>0$.
\end{cor}
\begin{proof}
  If $\spt f\cap B_{2s}(x)=\emptyset$ then $
  \inf_{z\in B_{s}(x)} Mf(z) \lesssim s^{-1}\|f\|_{L^{1}}$.
  Using this fact and interpolating between the bounds from
  Proposition \ref{prop:local-L2} and Lemma \ref{lem:locality-decay}
  we obtain the required inequality. 
\end{proof}

Fix $p\in(1,\infty]$ and $q\in\big(\max(p';2),\infty\big]$ and let
$\overline{p}\in\big(1,\min(p;2)\big)$ such that $q>\overline{p}'$.
We will now show that 
\begin{align*}\num\label{eq:local-lq-energy}
  \|F \1_{\mb X \setminus K_{f,\lambda,\overline{p}}}\|_{\L^{q}(S_{e})}\lesssim \lambda.
\end{align*}
Since $\nu(K_{f,\lambda,\overline{p}})\lesssim \lambda^{-p}\|f\|_{L^{p}}^{p}$ this
would prove \eqref{eq:e-bounds-w}.

Let us consider a strip
$D(x,s)\in \mb D$ and suppose that
$D(x,s)\not \subset K_{f,\lambda,\overline{p}}$, otherwise the estimate is
trivial.  We have $B_{5s}(x)\not\subset \mc K_{f,\lambda,\overline{p}}$. For an
$N>1$ large enough to be chosen later let us set
\begin{align*}
  f(x)=f_{0}(x)+\sum_{k=1}^{\infty} f_{k}(x) = f(x)\upsilon\Big(\frac{x-x_{0}}{5s}\Big)+\sum_{k=1}^{\infty}f(x)
  \gamma\left(\frac{x-x_{0}}{5s2^{Nk}}\right) 
\end{align*}
where $\gamma(\cdot)=\upsilon(\cdot/2^{N})-\upsilon(\cdot)$ with
\begin{align*}
  &  \upsilon \in C_{c}^{\infty}(B_{2})
  &&\upsilon\geq 0
  &&\upsilon=1 \text{
     on } B_{1}.
\end{align*}
Let $F_{k}$ be associated to $f_{k}$ via the embedding \eqref{eq:F} and
let $K_{f_{k},\lambda,\overline{p}}$ be as in
\eqref{eq:maximal-superlevel}.

Since $K_{f_{0},\lambda,\overline{p}}\subset K_{f,\lambda,\overline{p}}$
we have  that $\|f_{0}\|_{L^{\overline{p}}}s^{-1/\overline{p}}\lesssim
\lambda$ and 
\begin{align*}\num\label{eq:F0}
  \|F_{0}\,\1_{\mb X\setminus K_{f,\lambda,\overline{p}}} \,\1_{D(x,s)}\|_{L^{q}(S_{e})} \lesssim
  \lambda^{1-\overline{p}/q} \|f_{0}\|_{L^{\overline{p}}}^{\overline{p}/q}
  \lesssim \lambda  s^{1/q}
\end{align*}
by Theorem \ref{thm:thm-diou}

Since
$K_{f_{k},\lambda,\overline{p}}\subset K_{f,\lambda,\overline{p}}\not\supset B_{5s}(x)$ one
has $ \|f_{k}\|_{L^{\overline{p}}}\lesssim \lambda \nu(D(x,s))^{1/\overline{p}} 2^{Nk/\overline{p}}$ and by
Corollary \ref{cor:local} we have that
\begin{align*}\num\label{eq:Fk}
&  \|F_{k}\,\1_{\mb X\setminus K_{f,\lambda,p}}\,\1_{D(x,s)}\|_{L^{q}(S_{e})}s^{-1/q}\lesssim
  2^{-2Nk}2^{Nk/\overline{p}}\lambda\lesssim 2^{-Nk}\lambda. 
\end{align*}
By quasi-subadditivity we can add up \eqref{eq:F0} and \eqref{eq:Fk}
to obtain
\begin{align*}
  \frac{\|F\,\1_{\mb X\setminus
  K_{f,\lambda,p}}\,\1_{D(x,s)}\|_{L^{q}(S_{e})}}{\nu(D(x,s))^{1/q}}\lesssim \lambda.
\end{align*}
Since $D(x,s)$ is arbitrary this implies \eqref{eq:local-lq-energy}.

\subsection{Non-iterated bounds}

We conclude by explaining that for $r\in(2,\infty]$ and $p\in(2,r)$
and simpler embedding bounds on the maps $f\mapsto F$ and
$a\mapsto \mb A$ are sufficient to prove boundedness on $L^{p}(\R)$ of
the Variational Carleson Operator \eqref{eq:var-carleson-op} and thus
also \eqref{eq:carleson-op}.

Hereafter we work with the non-iterated outer measure space
$(\mb X,\mu)$.  The energy embedding map satisfies the $L^{p}$ bounds
\begin{align*}\num\label{eq:F-bounds-local-L2}
  &\|F\|_{L^{p}(S_{e})}\lesssim \|f\|_{L^{p}}&&p\in(2,\infty].
\end{align*}
This follows directly from Proposition \ref{prop:local-L2}  by taking
$s$ arbitrarily large.

Similarly, in Proposition \ref{prop:aux-map-non-iter} we have shown
that the auxiliary embedding satisfies
\begin{align*}\num\label{eq:aux-bounds-local-L2}
  &\|\mb M\|_{L^{p'}(S_{m})}\lesssim \|\a\|_{L^{p'}(l^{r'})}&&p'\in(r',\infty].
\end{align*}
and thus, by Proposition \ref{prop:size-control} we have that the
variational mass embedding also satisfies such bounds:
\begin{align*}\num\label{eq:A-bounds-local-L2}
  &\|\mb A\|_{L^{p'}(S_{m})}\lesssim \|\a\|_{L^{p'}(l^{r'})}&&p'\in(r',\infty].
\end{align*}

It follows by the outer Hölder inequality \ref{prop:outer-holder} that
\begin{align*}
\Big|  \iiint_{\mb X} F(y,\eta,t) \mb A(y,\eta,t)\dd y\dd \eta\dd
  t \Big| \lesssim \|F\|_{L^{p}(S_{e})}\|\mb A\|_{L^{p'}(S_{m})}.
\end{align*}
Using \eqref{eq:F-bounds-local-L2} and \eqref{eq:A-bounds-local-L2}
and the wave-packet domination \eqref{eq:var-duality-wavepacket} it
follows that \eqref{eq:var-carleson-op} is bounded on $L^{p}(\R)$.

In conclusion we remark that the iterated outer-measure $L^p$ spaces
that were introduced provide an effective way of capturing the spatial
locality property of the embedding maps. Both the proof of Theorem
\ref{thm:max-var-mass-bounds} and of Theorem
\ref{thm:energy-boundedness} rely on first obtaining non-iterated
bounds (see Propositions \ref{prop:aux-map-non-iter} and
\ref{prop:local-L2}) and then using a locality lemma (see Lemmata
\ref{lem:locality-M-Om-R} and
\ref{lem:locality-decay} ) and  a projection lemma (see
Lemma \ref{lem:mass-projection} and Lemma 7.8 of
\cite{di2015modulation}) to bootstrap the full result. 


\printbibliography

\end{document}